\documentclass[12pt,oneside,reqno]{amsart}
\usepackage{bbding}
\usepackage{txfonts}
\usepackage{bbm}
\usepackage{amsmath}
\usepackage{graphicx}
\usepackage{mathrsfs}
\usepackage{stmaryrd}
\usepackage{amsfonts}
\usepackage{enumerate,amsmath,amssymb}

\numberwithin{equation}{section}

\newcommand{\be}{\begin{eqnarray}}
\newcommand{\ee}{\end{eqnarray}}
\newcommand{\ce}{\begin{eqnarray*}}
\newcommand{\de}{\end{eqnarray*}}
\newtheorem{theorem}{Theorem}[section]
\newtheorem{lemma}[theorem]{Lemma}
\newtheorem{remark}[theorem]{Remark}
\newtheorem{definition}[theorem]{Definition}
\newtheorem{proposition}[theorem]{Proposition}
\newtheorem{Examples}[theorem]{Example}
\newtheorem{corollary}[theorem]{Corollary}

\def\eps{\varepsilon}

\def\ind{r}

\def\osc{\mathrm{osc}}
\def\med{\mathrm{med}}

\def\p{\partial}

\def\[{{\Big[}}
\def\]{{\Big]}}
\def\<{{\langle}}
\def\>{{\rangle}}
\def\({{\Big(}}
\def\){{\Big)}}

\def\bx{{\mathbf{x}}}

\def\dif{{\mathord{{\rm d}}}}

\def\min{{\mathord{{\rm min}}}}

\def\no{\nonumber}
\def\={&\!\!=\!\!&}
\def\bt{\begin{theorem}}
\def\et{\end{theorem}}
\def\bl{\begin{lemma}}
\def\el{\end{lemma}}
\def\br{\begin{remark}}
\def\er{\end{remark}}

\def\bd{\begin{definition}}
\def\ed{\end{definition}}
\def\bp{\begin{proposition}}
\def\ep{\end{proposition}}
\def\bc{\begin{corollary}}
\def\ec{\end{corollary}}
\def\bx{\begin{Examples}}
\def\ex{\end{Examples}}

\def\cF{{\mathcal F}}

\def\cH{{\mathcal H}}

\def\cL{{\mathcal L}}
\def\cM{{\mathcal M}}

\def\cP{{\mathcal P}}

\def\cS{{\mathcal S}}
\def\cT{{\mathcal T}}

\def\mD{{\mathbb D}}
\def\mE{{\mathbb E}}

\def\mH{{\mathbb H}}

\def\mM{{\mathbb M}}
\def\mN{{\mathbb N}}

\def\mQ{{\mathbb Q}}
\def\mR{{\mathbb R}}
\def\mS{{\mathbb S}}

\def\mW{{\mathbb W}}

\def\sB{{\mathscr B}}
\def\sC{{\mathscr C}}
\def\sD{{\mathscr D}}

\def\sF{{\mathscr F}}

\def\sL{{\mathscr L}}

\def\sT{{\mathscr T}}

\def\geq{\geqslant}
\def\leq{\leqslant}

\allowdisplaybreaks

\begin{document}

\title{$L^p$-maximal regularity of nonlocal parabolic equation and applications$^*$}

\date{}
\author{Xicheng Zhang}

\thanks{{\it Keywords: }$L^p$-regularity, L\'evy process, Krylov's estimate, sharp function, critical Burger's equation}

\thanks{$*$This work is supported by NSFs of China (No. 10971076).}

\address{
School of Mathematics and Statistics,
Wuhan University, Wuhan, Hubei 430072, P.R.China,\\
Email: XichengZhang@gmail.com
 }

\begin{abstract}
By using Fourier's transform and Fefferman-Stein's theorem, we investigate the $L^p$-maximal regularity
of nonlocal parabolic and elliptic equations with singular and non-symmetric L\'evy operators,
and obtain the unique strong solvability of the corresponding nonlocal parabolic and elliptic equations, where
the probabilistic representation plays an important role. As a consequence,
a characterization for the domain of pseudo-differential operators of L\'evy type with singular kernels
is given in terms of the Bessel potential spaces. As a byproduct, we also show that a large class of non-symmetric L\'evy operators
generates an analytic semigroup in $L^p$-space. Moreover, as applications, we prove a Krylov's estimate
for stochastic differential equation driven by Cauchy processes (i.e. critical diffusion processes),
and also obtain the global well-posedness to a class of quasi-linear first order parabolic
system with critical diffusion. In particular, critical Hamilton-Jacobi equation and multidimensional
critical Burger's equation are uniquely solvable and the smooth solutions are obtained.
\end{abstract}

\maketitle \rm

\section{Introduction}

Consider the following Cauchy problem of fractional Laplacian heat equation
in the domain $[0,\infty)\times\mR^d$ with $\alpha\in(0,2)$ and $\lambda\geq 0$:
\begin{align}
\p_t u+(-\Delta)^{\frac{\alpha}{2}}u+b\cdot\nabla u+\lambda u=f,\ u(0)=\varphi, \label{Eq1}
\end{align}
where $b:[0,\infty)\times\mR^d\to\mR^d$ is a measurable vector field, $f:[0,\infty)\times\mR^d\to\mR$ and $\varphi:\mR^d\to\mR$
are two measurable functions, and $(-\Delta)^{\frac{\alpha}{2}}$ is the fractional Laplacian
(also called L\'evy operator)  defined  by
\begin{align}
(-\Delta)^{\frac{\alpha}{2}}u=\cF^{-1}(|\cdot|^{\alpha}\cF(u)),\ \ u\in\cS(\mR^d),\label{ET7}
\end{align}
where $\cF$ (resp. $\cF^{-1}$) denotes the Fourier (resp. inverse) transform, $\cS(\mR^d)$ is the Schwartz class of smooth
real or complex-valued rapidly decreasing functions.

Let $(L_t)_{t\leq 0}$ be a symmetric and rotationally invariant
$\alpha$-stable process. Let $b,f\in C^\infty_b([0,\infty)\times\mR^d)$ and $X_{t,s}(x)$
solve the following stochastic differential equation (SDE):
$$
X_{t,s}(x)=x+\int^s_tb(-r,X_{t,r}(x))\dif r+\int^s_t\dif L_r,\ \ t\leq s\leq 0,\ x\in\mR^d.
$$
It is well-known that for $\varphi\in C^\infty_b(\mR^d)$,
the unique solution of equation (\ref{Eq1}) can be represented by Feyman-Kac formula as
(see Theorem \ref{EY2} below):
\begin{align}
u(t,x)=\mE\varphi(X_{-t,0}(x))+\mE\left(\int^0_{-t}e^{-\lambda(s+t)} f(-s, X_{-t,s}(x))\dif s\right),\ \
t\geq 0.\label{Prob}
\end{align}
In connection with this representation, the first order term $b\cdot\nabla u$
is also called the drift term, and the fractional
Laplacian term $(-\Delta)^{\frac{\alpha}{2}}u$ is also called the diffusion term.

Let now $u(t,x)$ satisfy (\ref{Eq1}). For $\ind>0$ and $(t,x)\in[0,\infty)\times\mR^d$, define
$$
u^\ind(t,x):=\ind^{-\alpha} u(\ind^\alpha t, \ind x),\ \
b^\ind(t,x):= b(\ind^\alpha t,\ind x),\ \
f^\ind(t,x):= f(\ind^\alpha t, \ind x),
$$
then it is easy to see that $u^\ind$ satisfies
\begin{align}
\p_t u^\ind+(-\Delta)^{\frac{\alpha}{2}}u^\ind+\ind^{\alpha-1}(b^\ind\cdot\nabla u^\ind)+\lambda r^\alpha u^r=f^\ind.
\label{Sc1}
\end{align}
If one lets $\ind\to 0$, this scaling property leads to the following classification:
\begin{enumerate}[$\bullet$]
\item (Subcritical case: $\alpha\in(1,2)$) The drift term is controlled by the diffusion
term at small scales.

\item (Critical case: $\alpha=1$) The fractional Laplacian has the same order as the first order term.

\item (Supercritical case: $\alpha\in(0,1)$) The effect of the drift term is stronger
than the diffusion term at small scales.
\end{enumerate}

In recent years there are great interests to study the above nonlocal equation, since
it has appeared in numerous disciplines, such as quasi-geostrophic fluid dynamics (cf. \cite{Con, Ca-Va}),
stochastic control problems (cf. \cite{So}), nonlinear filtering with jump (cf. \cite{Mi-Pr2}),
mathematical finance (cf. \cite{Be-Ka-Re}), anomalous diffusion in semiconductor growth (cf. \cite{Wo}), etc.
In \cite{Dr-Im}, Droniou and Imbert studied the first order Hamilton-Jacobi equation
with fractional diffusion $(-\Delta)^{\frac{\alpha}{2}}$ basing upon a ``reverse maximal principle''. Therein,
when $\alpha\in(1,2)$, the classical solution was obtained; when $\alpha\in(0,2)$,
the existence and uniqueness of viscosity solutions in the class of Lipschitz functions were also established.
In \cite{Ca-Va}, Caffarelli and Vasseur established the global well-posedness of critical dissipative quasi-geostrophic
equation (see also \cite{Ki-Na-Vo} for a simple proof in the periodic and two dimensional case). On the other hand,
H\"older regularity theory for
the viscosity solutions of fully non-linear and nonlocal elliptic equations are also developed by Caffarelli and Silvestre \cite{Ca-Si},
and Barles, Chasseigne and Imbert \cite{Ba-Ch-Im}, see also the series of works of Silverstre \cite{Si0,Si3,Si1,Si2}, etc.
We emphasize that the arguments in \cite{Ca-Si} and \cite{Ba-Ch-Im} are different: the former is
based on the Alexandorff-Backelman-Pucci's (ABP) estimate, and the latter is based on the Ishii-Lions' simple method.

The purpose of this paper is an attempt to develop an $L^p$-regularity theory
for nonlocal equations with general L\'evy operators. We describe it as follows.
Let $\nu$ be a L\'evy measure in $\mR^d$, i.e., a $\sigma$-finite measure satisfying $\nu(\{0\})=0$ and
$$
\int_{\mR^d}\min(1,|y|^2)\nu(\dif y)<+\infty.
$$
For $\alpha\in(0,2)$, we write
$$
y^{(\alpha)}:=1_{\alpha\in(1,2)}y+1_{\alpha=1}y1_{|y|\leq 1}.
$$
In this article we are mainly concerned with the following pseudo-differential operator of L\'evy type:
\begin{align}
\cL^{\nu} f(x):=\int_{\mR^d}[f(x+y)-f(x)-
y^{(\alpha)}\cdot\nabla f(x)]\nu(\dif y),\ \ f\in\cS(\mR^d),\label{Laa}
\end{align}
where $\nu$ satisfies
\begin{align}
\nu^{(\alpha)}_1(B)\leq \nu(B)\leq\nu^{(\alpha)}_2(B),\ \ B\in \sB(\mR^d),\label{EL2}
\end{align}
and
\begin{align}
1_{\alpha=1}\int_{r\leq |y|\leq R}y\nu(\dif y)=0, \ \ 0<r<R<+\infty.\label{EL22}
\end{align}
Here, $\nu^{(\alpha)}_i, i=1,2$ are the L\'evy measures of two $\alpha$-stable processes taking the form
\begin{align}
\nu^{(\alpha)}_i(B):=\int_{\mS^{d-1}}\left(\int^\infty_0
\frac{1_B(r\theta)\dif r}{r^{1+\alpha}}\right)\Sigma_i(\dif\theta),\label{Eq4}
\end{align}
where $\mS^{d-1}=\{\theta\in\mR^d: |\theta|=1\}$ is the unit sphere in $\mR^d$, and $\Sigma_i$
called the spherical part of $\nu^{(\alpha)}_i$ is a finite measure on $\mS^{d-1}$.
We remark that condition (\ref{EL22}) is a common assumption in the critical case (see \cite{Mi-Pr1, Do-Ki}),
and is clearly satisfied when $\nu$ is symmetric.

One of the aims of the present paper is to determine $\sD^p(\cL^\nu)$, the domain of the L\'evy operator $\cL^\nu$ in $L^p$-space.
We shall prove that under (\ref{EL2}) and (\ref{EL22}), if $\nu^{(\alpha)}_1$ is nondegenerate
(see Definition \ref{Def} below), then for any $p\in(1,\infty)$,
$$
\sD^p(\cL^\nu)=\mH^{\alpha,p},
$$
where $\mH^{\alpha,p}$ is the $\alpha$-order Bessel potential space.
When $\nu(\dif y)=a(y)\dif y/|y|^{d+\alpha}$
with $c_1\leq |a(y)|\leq c_2$, this characterization was obtained recently by Dong and Kim \cite{Do-Ki}.
It is remarked that the technique in \cite{Ba-Ch-Im} was used by Dong and Kim to derive
some local H\"older estimate for nonlocal elliptic equation in order to prove their characterization.
However, the following sum of nonlocal operators is not covered by \cite{Do-Ki}:
$$
\cL f(x)=\sum_{i=1}^d\int_{\mR}\frac{f(x_1,\cdots,x_{i-1},
x_i+y_i,x_{i+1},\cdots x_d)-f(x)-y_i^{(\alpha)}\cdot\p_i f(x)}{|y_i|^{1+\alpha}}\dif y_i,
$$
since in this case, the L\'evy measure (or the L\'evy symbol) is very singular
(or non-smooth) (see Remark \ref{Re2}).
Notice that if the L\'evy symbol is smooth and its derivatives satisfy suitable conditions, the above characterization
falls into the classical multiplier theorems about pseudo-differential operators (cf. \cite{St0, Ja}).
We also mention that Farkas, Jacob and Schilling \cite[Theorem 2.1.15]{Fa-Ja-Sc}
gave another characterization for $\sD^p(\cL^\nu)$ in terms of the so called $\psi$-Bessel potential space,
where $\psi$ is the symbol of $\cL^\nu$.

The strategy for proving the above characterization is to prove the
following Littlewood-Paley type inequality: for any $p\in(1,\infty)$, there exists a $C>0$ such that for any $\lambda\geq 0$,
$f\in L^p(\mR^+\times\mR^d)$,
$$
\int^\infty_0\left\|\cL^{\nu_2}\int^t_0e^{-\lambda(t-s)}\cP^{\nu_1}_{t-s}f(s,\cdot)\dif s\right\|_p^p\dif t
\leq C\int^\infty_0\|f(t,\cdot)\|^p_p\dif t,
$$
where $\nu_1,\nu_2$ are two L\'evy measures satisfying (\ref{EL2}) and (\ref{EL22}), and
$(\cP^{\nu_1}_t)_{t\geq 0}$ is the semigroup associated with $\cL^{\nu_1}$.
Indeed, this estimate is the key ingredient in $L^p$-theory of PDE (see \cite{La-So-Ur, Kr3}),
and corresponds to the optimal regularity of nonlocal parabolic equation. Likewise \cite{Do-Ki}, when
$\nu(\dif y)=a(y)\dif y/|y|^{d+\alpha}$ with smooth and $0$-homogeneous $a(y)$ and
$c_1\leq |a(y)|\leq c_2$, Mikulevicius and Pragarauskas \cite{Mi-Pr1}
proved this type of estimate by showing some weak $(1,1)$-type estimate. In a different way,
the proof given here is based on Fourier's transform and Fefferman-Stein's theorem about sharp functions (cf. \cite{Kr1,Kr3}).
We stress that probabilistic representation
(\ref{Prob}) will play an important role in reducing the general nonhomogeneous operator to homogeneous operator
(see Step 1 in the proof of Theorem \ref{Th1}).

Another aim of this paper is to solve the linear and quasi-linear first order nonlocal
parabolic equation with critical diffusions
in the $L^p$-sense rather than the viscosity sense \cite{Dr-Im}.
The critical case is specially interesting not only because it appears naturally in quasi-geostrophic flows,
but also it is an attractive object in mathematics. In particular, we care about the
following multidimensional critical Burger's equation:
\begin{align}
\p_t u+(-\Delta)^{\frac{1}{2}}u+u\cdot \nabla u=0,\ \ u(0)=\varphi.\label{Eq0}
\end{align}
In one dimensional case, this equation has a natural variational formulation and admits a unique global smooth solution
(see \cite{Bi-Fu-Wo,Ki-Na-Sh}) under some regularity assumption on $\varphi$.
In multidimensional case, the local well-posedness of Burger's equation is relatively easy (cf. \cite{Ja-Po-Wu, Zh4}).
However, the global well-posedness of equation (\ref{Eq0}) is a challenge problem.
The reason lies in two aspects: on one hand, there is no energy inequality and thus, the variational method
seems not to be applicable; on the other hand, the first order term has the same order as the diffusion term.
In fact, Kiselev, Nazarov and Schterenberg \cite{Ki-Na-Sh} have showed the existence of blow up solutions
for $1$-D supercritical Burger's equation.
The idea here is to establish some apriori H\"older estimate for equation (\ref{Eq1})
and then use the classical method of freezing coefficients. In \cite{Si2},
Silvestre proved an apriori H\"older estimate for equation (\ref{Eq1}) with only bounded measurable $b$.
This is the key point for us. However, the assumption of scale invariance on L\'evy operators
seems to be crucial in \cite{Si2} since the proof is by the iteration of the diminish of oscillation at all scales.
As above, we shall use probabilistic representation (\ref{Prob})
like a perturbation argument to extend Silvestre's estimate to
the more general non-homogeneous L\'evy operator (see Corollary \ref{Cor4}).

This paper is organized as follows. In Section 2, we prepare some lemmas and recall some facts for later use.
In Section 3, the basic maximum principles for nonlocal parabolic and elliptic equation are proved.
In Section 4, we prove the main Theorem \ref{Th1}, and give a comparison result between two L\'evy operators.
In particular, we show that $(\cP^\nu_t)_{t\geq 0}$ forms an analytic semigroup in $L^p$-space.
In Section 5, we prove the existence of a unique strong solution for
the first order nonlocal parabolic equation with critical diffusion and various coefficients.
As an application, we also prove a Krylov's estimate for critical diffusion processes.
We mention that in one dimensional and subcritical case, such type of estimate was firstly proved by Kurenok \cite{Ku};
and in multidimensional and subcritical case, it was proved in \cite{Zh3}.
In Section 6, we investigate quasi-linear first order nonlocal parabolic system, and get the existence of
smooth solutions and strong solutions. In particular, the global solvability
of equation (\ref{Eq0}) is obtained.

\vspace{5mm}

{\bf Notations:} We collect some frequently used notations below for the reader's convenience.

\begin{enumerate}[$\bullet$]
\item $\mR^+:=(0,\infty)$, $\mR^+_0:=[0,\infty)$. For a complex number $z$, Re($z$) (Im($z$)): real (image) part of $z$.
\item $\cS(\mR^d)$: the Schwartz class of smooth real or complex-valued rapidly decreasing functions.
$C^\infty_b(\mR^d)$ (resp. $C^k_b(\mR^d)$, $C^\infty_0(\mR^d)$):
the space of all bounded smooth functions with bounded derivatives of all orders
(resp. up to $k$-order, with compact support).
\item $\cF$ and $\cF^{-1}$: Fourier's transform and Fourier's inverse transform.
\item $\nu$: L\'evy measure; $\nu^{(\alpha)}$: the L\'evy measure of $\alpha$-stable process; $\Sigma$: a finite measure
on $\mS^{d-1}$, called the spherical part of $\nu^{(\alpha)}$.
\item $L^\nu_t$: the L\'evy process associated with L\'evy measure $\nu$; $\cP^\mu_t$: the semigroup associated with $L^\mu_t$.
$\cL^\nu$: the generator of $L^\mu_t$, $\cL^{\nu*}$: the adjoint operator of $\cL^\nu$; $p^\nu_t$: the heat kernel of $\cL^{\nu*}$.
\item $B_r(x_0):=\{x:\in\mR^d: |x-x_0|\leq r\}$,\ \ $B_r:=B_r(0)$,\ \ $B^c_r$: the complement of $B_r$.
\item $\mH^{\alpha,p}$: Bessel potential space; $\mW^{\alpha,p}$: Sobolev-Slobodeckij space; $\mW^\infty:=\cap_{k,p}\mW^{k,p}$.
\item $\omega_b$: the continuous modulus function of $b$, i.e., $\omega_b(s):=\sup_{|x-y|\leq s}|b(x)-b(y)|$.
\item $\cH^\beta$: the space of H\"older continuous functions with the norm $\sum_{k=0}^{[\beta]}\|\nabla^k f\|_\infty+\|\nabla^{[\beta]}f\|_{\cH^\beta}$,
where $[\beta]$ denotes the integer part of $\beta$,
and $\|\nabla^{[\beta]}f\|_{\cH^\beta}:=\sup_{|x-y|\leq 1}|\nabla^{[\beta]}f(x)-\nabla^{[\beta]}f(y)|/|x-y|^\beta$.
\item $(\rho_\eps)_{\eps\in(0,1)}$: a family of mollifiers in $\mR^d$ with $\rho_\eps(x)=\eps^{-d}\rho(\eps^{-1}x)$,
where $\rho$ is a nonnegative smooth function with support in $B_1$ and satisfies $\int_{\mR^d}\rho(x)\dif x=1$.
\end{enumerate}

{\bf Convention:} The letter $C$ with or without subscripts will denote an unimportant constant.
The inner product in Euclidean space is denoted by ``$\cdot$''.

\section{Preliminaries}

For $\alpha\in(0,2)$, let $\nu$ be a L\'evy measure in $\mR^d$ and satisfy (\ref{EL2}) and (\ref{EL22}).
Let $(L^\nu_t)_{t\geq 0}$ be the $d$-dimensional L\'evy process,
a stationary and independent increment process defined on some probability space $(\Omega,\sF,P)$,
with characteristic function
\begin{align}
\mE e^{\mathrm{i}\xi\cdot L^{\nu}_t}=e^{-t\psi_\nu(\xi)},\ \ \xi\in\mR^d,\label{Fr}
\end{align}
where $\psi_v$ is the L\'evy exponent with the form by L\'evy-Khintchine's formula (cf. \cite{Ap,Sa}),
\begin{align}
\psi_\nu(\xi):=\int_{\mR^d}(1+\mathrm{i}
\xi\cdot y^{(\alpha)}-e^{\mathrm{i}\xi\cdot y})\nu(\dif y).\label{Fr1}
\end{align}
Let $\nu^{(\alpha)}$ take the form (\ref{Eq4}) and satisfy (\ref{EL22}). It is well-known that
$(L^{\nu^{(\alpha)}}_t)_{t\geq 0}$ is a $d$-dimensional $\alpha$-stable process and has
the following self-similarity (cf. \cite[Proposition 13.5 and Theorem 14.7]{Sa}):
\begin{align}
(L^{\nu^{(\alpha)}}_{\ind t})_{t\geq 0}
\stackrel{(d)}{=}(\ind^{1/\alpha} L^{\nu^{(\alpha)}}_t)_{t\geq 0},\ \ \forall \ind>0,\label{Sc}
\end{align}
where $\stackrel{(d)}{=}$ means that the two processes have the same laws. Moreover, from expression (\ref{Eq4}),
it is easy to see that for any $\beta\in(0,\alpha)$,
\begin{align}
\int_{\mR^d}\min(|y|^\beta,|y|^2)\nu^{(\alpha)}(\dif y)<+\infty,\label{Es6}
\end{align}
and
\begin{align}
\mathrm{Re}(\psi_{\nu^{(\alpha)}}(\xi))=\left(\int^\infty_0\frac{(1-\cos r)\dif r}{r^{1+\alpha}}\right)
\int_{\mS^{d-1}}|\xi\cdot\theta|^\alpha\Sigma(\dif\theta).\label{Es2}
\end{align}

The Feller semigroup associated with $(L^\nu_t)_{t\geq 0}$ is defined by
$$
\cP^{\nu}_t f(x):=\mE f(L^{\nu}_t+x),\ \ f\in\cS(\mR^d).
$$
The generator of $(\cP^{\nu}_t)_{t\geq 0}$ is then given by (cf. \cite[Theorem 3.3.3]{Ap})
\begin{align}
\cL^{\nu} f(x)=\int_{\mR^d}[f(x+y)-f(x)-y^{(\alpha)}\cdot\nabla f(x)]\nu(\dif y),\label{La}
\end{align}
i.e.,
\begin{align}
\p_t \cP^{\nu}_tf(x)=\cL^{\nu} \cP^{\nu}_tf(x)= \cP^{\nu}_t\cL^{\nu}f(x),\ \ t>0.\label{Ep4}
\end{align}
Moreover,
$$
\cF(\cL^{\nu} f)(\xi)=-\psi_\nu(\xi)\cdot\cF(f)(\xi),
$$
and $\psi_\nu$ is also called the L\'evy symbol of the operator $\cL^{\nu}$. From (\ref{Es2}), one sees
that if the spherical part $\Sigma$ of $\nu^{(\alpha)}$ is the uniform distribution
(equivalently, rotationally invariant) on $\mS^{d-1}$,
then $\psi_{\nu^{(\alpha)}}(\xi)=c_{d,\alpha}|\xi|^\alpha$ for some constant $c_{d,\alpha}>0$, and thus,  by (\ref{ET7}),
\begin{align}
-\cL^{\nu^{(\alpha)}} f(x)=c_{d,\alpha}(-\Delta)^{\frac{\alpha}{2}}f(x).\label{EW5}
\end{align}
On the other hand, from expression (\ref{La}) and assumption (\ref{EL22}), it is easy to see that $\cL^{\nu}$ has the following invariance:
\begin{enumerate}[$\bullet$]
\item For $z\in\mR^d$, define $f_z(x):=f(z+x)$, then
\begin{align}
\cL^{\nu} f_z(x)=\cL^{\nu} f_x(z),\ \
\|\cL^{\nu} f_z\|_p=\|\cL^{\nu} f\|_p,\label{EU1}
\end{align}
where $p\geq 1$ and $\|\cdot\|_p$ denotes the usual $L^p$-norm in $\mR^d$.
\item For $\ind>0$, define $f_\ind(x):=f(\ind x)$, then
\begin{align}
\cL^{\nu} f(\ind x)=\cL^{\nu(r\cdot)} f_\ind(x)=\ind^{-\alpha}\cL^{r^\alpha\nu(r\cdot)} f_\ind(x).\label{EU2}
\end{align}
We remark that $r^\alpha\nu^{(\alpha)}(r\cdot)=\nu^{(\alpha)}$ by (\ref{Eq4}).
\item $\cL^{\nu}(C^\infty_b(\mR^d))\subset C^\infty_b(\mR^d)$, and for any $k\geq 2$,
$\cL^\nu: C^k_b(\mR^d)\to C^{k-2}_b(\mR^d)$ is a continuous linear operator,
where $C^\infty_b(\mR^d)$ (resp. $C^k_b(\mR^d)$) is the space of all bounded smooth
functions with bounded derivatives of all orders (resp. up to $k$-order).
\end{enumerate}
The adjoint operator of $\cL^{\nu}$ is given by
\begin{align}
\cL^{\nu*} f(x)=\int_{\mR^d}[f(x-y)-f(x)+
y^{(\alpha)}\cdot\nabla f(x)]\nu(\dif y),\label{Ad}
\end{align}
i.e.,
$$
\int_{\mR^d}\cL^{\nu} f(x)\cdot g(x)\dif x=\int_{\mR^d}f(x)\cdot\cL^{\nu*} g(x)\dif x,\ \ f,g\in\cS(\mR^d).
$$
Clearly, $\cL^{\nu*}=\cL^{\nu(-)}$, where $\nu(-)$ denotes the L\'evy measure $\nu(-\dif y)$.
\bd
Let $\nu_1$ and $\nu_2$ be two Borel measures. We say that $\nu_1$ is less than $\nu_2$ if
$$
\nu_1(B)\leq \nu_2(B),\ \ B\in\sB(\mR^d),
$$
and we simply write $\nu_1\leq \nu_2$ in this case.
\ed
\bl
Let $\nu$ be a L\'evy measure less than $\nu^{(\alpha)}$ for some $\alpha\in(0,2)$,
where $\nu^{(\alpha)}$ takes the form (\ref{Eq4}). We also assume (\ref{EL22}) for $\nu$.
Then for some $\kappa_0>0$,
\begin{align}
|\psi_{\nu}(\xi)|\leq \kappa_0|\xi|^\alpha,\ \ \xi\in\mR^d.\label{Es22}
\end{align}
\el
\begin{proof}
Write $\hat\xi:=\xi/|\xi|$. For $\alpha\in(1,2)$, by the definitions of $\psi_\nu$ and $\nu^{(\alpha)}$, we have
\begin{align*}
|\mathrm{Im}(\psi_\nu(\xi))|&\stackrel{(\ref{Fr1})}{\leq}\int_{\mR^d}|\xi\cdot y-\sin(\xi\cdot y)|\nu(\dif y)
\leq\int_{\mR^d}|\xi\cdot y-\sin(\xi\cdot y)|\nu^{(\alpha)}(\dif y)\\
&\stackrel{(\ref{Eq4})}{=}\int_{\mS^{d-1}}\int^\infty_0\frac{|\xi\cdot (r\theta)
-\sin(\xi\cdot r\theta)|}{r^{1+\alpha}}\dif r\Sigma(\dif\theta)\\
&=|\xi|^\alpha\int_{\mS^{d-1}}\int^\infty_0\frac{|\hat\xi\cdot r\theta
-\sin(\hat\xi\cdot r\theta)|}{r^{1+\alpha}}\dif r\Sigma(\dif\theta)\leq C|\xi|^\alpha.
\end{align*}
For $\alpha=1$, by (\ref{EL22}), we have
\begin{align*}
|\mathrm{Im}(\psi_\nu(\xi))|&=\left|\int_{\mR^d}(\xi\cdot y1_{|y|\leq|\xi|^{-1}}-\sin(\xi\cdot y))\nu(\dif y)\right|\\
&\leq\int_{\mR^d}|\xi\cdot y1_{|y|\leq|\xi|^{-1}}-\sin(\xi\cdot y)|\nu^{(1)}(\dif y)\\
&=\int_{\mS^{d-1}}\int^\infty_0\frac{|\xi\cdot (r\theta)1_{r\leq|\xi|^{-1}}
-\sin(\xi\cdot r\theta)|}{r^2}\dif r\Sigma(\dif\theta)\\
&=|\xi|\int_{\mS^{d-1}}\int^\infty_0\frac{|\hat\xi\cdot r\theta1_{r\leq 1}
-\sin(\hat\xi\cdot r\theta)|}{r^2}\dif r\Sigma(\dif\theta)\leq C|\xi|.
\end{align*}
For $\alpha\in(0,1)$, we have
\begin{align*}
|\mathrm{Im}(\psi_\nu(\xi))|&\leq\int_{\mR^d}|\sin(\xi\cdot y)|\nu(\dif y)
\leq\int_{\mR^d}|\sin(\xi\cdot y)|\nu^{(\alpha)}(\dif y)\\
&=|\xi|^\alpha\int_{\mS^{d-1}}\int^\infty_0\frac{|\sin(\hat\xi\cdot r\theta)|}
{r^{1+\alpha}}\dif r\Sigma(\dif\theta)\leq C|\xi|^\alpha.
\end{align*}
Thus, combining with (\ref{Es2}), we obtain (\ref{Es22}).
\end{proof}

For $k\in\mN$ and $p\in[1,\infty]$, let $\mW^{k,p}$ be the usual Sobolev space with the norm
$$
\|f\|_{k,p}:=\sum_{j=0}^k\|\nabla^j f\|_p,
$$
where $\nabla^j$ denotes the $j$-order gradient.

We need the following simple interpolation result.
\bl\label{Le4}
Let $p\in[1,\infty]$ and $\beta\in[0,1]$. For any $f\in\mW^{1,p}$ and $y\in\mR^d$, we have
\begin{align}
\|f(\cdot+y)-f(\cdot)\|_p\leq(2\|f\|_p)^{1-\beta}(\|\nabla f\|_p|y|)^\beta.\label{EW44}
\end{align}
\el
\begin{proof}
Observing that for $f\in\cS(\mR^d)$,
$$
|f(x+y)-f(x)|\leq |y|\int^1_0 |\nabla f|(x+sy)\dif s,
$$
by a density argument, we have for any $f\in\mW^{1,p}$,
$$
\|f(\cdot+y)-f(\cdot)\|_p\leq \|\nabla f\|_p|y|.
$$
Thus, for any $\beta\in[0,1]$,
$$
\|f(\cdot+y)-f(\cdot)\|_p\leq (2\|f\|_p)\wedge(\|\nabla f\|_p|y|)\leq (2\|f\|_p)^{1-\beta}(\|\nabla f\|_p|y|)^\beta.
$$
The result follows.
\end{proof}
The following lemma will be used to derive some asymptotic estimate of large time for the heat kernel of L\'evy operator
(see Corollary \ref{Cor1} below).
\bl\label{Le6}
Assume that L\'evy measure  $\nu$ is less than $\nu^{(\alpha)}$ for some $\alpha\in(0,2)$,
where $\nu^{(\alpha)}$ takes the form (\ref{Eq4}).
Then for any $p\in[1,\infty]$ and $f\in \mW^{2,p}$, we have
$$
\|\cL^{\nu} f\|_p\leq C\left\{
\begin{array}{llll}
 \|\nabla f\|_p^{1-\gamma}\|\nabla^2 f\|_p^\gamma+\|\nabla f\|_p^{1-\beta}\|\nabla^2 f\|^\beta_p,
&\alpha\in(1,2),&\gamma\in(\alpha-1,1], &\beta\in[0,\alpha-1), \\
\|\nabla f\|_p^{1-\gamma}\|\nabla^2 f\|^\gamma_p
+\|f\|_p^{1-\beta}\|\nabla f\|^\beta_p,& \alpha=1,& \gamma\in(0,1],&\beta\in[0,1),\\
\|f\|_p^{1-\gamma}\|\nabla f\|^\gamma_p+\|f\|_p^{1-\beta}\|\nabla f\|^\beta_p,&
\alpha\in(0,1),&\gamma\in(\alpha,1], &\beta\in[0,\alpha),
\end{array}
\right.
$$
where the constant $C$ depends only on $\alpha,\beta,\gamma$ and the L\'evy measure $\nu^{(\alpha)}$.
\el
\begin{proof}
Let us first look at the case of $\alpha\in(1,2)$. In this case, we have
$$
\cL^{\nu} f(x)=\int_{\mR^d}y\cdot\left(\int^1_0[\nabla f(x+sy)-\nabla f(x)]\dif s\right)\nu(\dif y)
$$
Since $\nu$ is bounded by $\nu^{(\alpha)}$,
by Minkowski's inequality and Lemma \ref{Le4}, we have for $\gamma\in(\alpha-1,1]$ and $\beta\in[0,\alpha-1)$,
\begin{align*}
\|\cL^{\nu} f\|_p&\leq (2\|\nabla f\|_p)^{1-\gamma}\|\nabla^2 f\|_p^\gamma
\int_{|y|\leq 1}|y|^{1+\gamma}\nu^{(\alpha)}(\dif y)
+(2\|\nabla f\|_p)^{1-\beta}\|\nabla^2 f\|^\beta_p
\int_{|y|>1}|y|^{1+\beta}\nu^{(\alpha)}(\dif y).
\end{align*}
In the case of $\alpha=1$, we similarly have for $\gamma\in(0,1]$ and  $\beta\in[0,1)$,
\begin{align*}
\|\cL^{\nu} f\|_p&\leq (2\|\nabla f\|_p)^{1-\gamma}\|\nabla^2 f\|^\gamma_p
\int_{|y|\leq1}|y|^{1+\gamma}\nu^{(1)}(\dif y)
+(2\|f\|_p)^{1-\beta}\|\nabla f\|^\beta_p
\int_{|y|>1}|y|^\beta\nu^{(1)}(\dif y).
\end{align*}
In the case of $\alpha\in(0,1)$, we have for $\gamma\in(\alpha,1]$ and $\beta\in[0,\alpha)$,
$$
\|\cL^{\nu} f\|_p\leq (2\|f\|_p)^{1-\gamma}\|\nabla f\|^\gamma_p
\int_{|y|\leq 1}|y|^{\gamma}\nu^{(\alpha)}(\dif y)+(2\|f\|_p)^{1-\beta}\|\nabla f\|^\beta_p
\int_{|y|>1}|y|^\beta\nu^{(\alpha)}(\dif y).
$$
The proof is complete by (\ref{Es6}).
\end{proof}
We also need the following estimate, which will be used frequently in localizing the nonlocal equation.
\bl\label{Le1}
Assume that L\'evy measure  $\nu$ is less than $\nu^{(\alpha)}$ for some $\alpha\in(0,2)$,
where $\nu^{(\alpha)}$ takes the form (\ref{Eq4}).
Let $\zeta\in \cS(\mR^d)$ and set $\zeta_z(x):=\zeta(x-z)$ for $z\in\mR^d$.

\begin{enumerate}[(i)]
\item For any $\beta\in(0\vee(\alpha-1),1)$ and $p\in[1,\infty)$, there exists a constant
$C=C(\nu^{(\alpha)},\beta,p,d)>0$ such that for all $ f\in\mW^{1,p}$,
\begin{align}
\left(\int_{\mR^d}\|\cL^{\nu}( f\zeta_z)
-(\cL^{\nu}f)\zeta_z\|_p^p\dif z\right)^{1/p}
\leq C\|\zeta\|_{2,p}\|f\|_p^{1-\beta}\|f\|_{1,p}^{\beta}.\label{EE3}
\end{align}
\item For any $\beta\in(0\vee(\alpha-1),1)$ and $\gamma\in[0,\alpha)$, there exists a constant
$C=C(\nu^{(\alpha)},\beta,\gamma,d)>0$ such that for any $p\in[1,\infty]$ and $f\in\cH^\beta$,
\begin{align}
\|\cL^{\nu}(f\zeta)-(\cL^{\nu}f)\zeta\|_p\leq
C\Big((\|\cL^\nu\zeta\|_p+\|\zeta\|^{1-\gamma}_p\|\nabla\zeta\|^\gamma_p)\|f\|_\infty
+\|\nabla\zeta\|_p\|f\|_{\cH^\beta}\Big),\label{EE4}
\end{align}
where $\|f\|_{\cH^\beta}:=\sup_{x\not=y, |x-y|\leq 1}|f(x)-f(y)|/|x-y|^\beta$, and for any
$p\in[1,\infty]$ and $f\in\mW^{1,p}$,
\begin{align}
\|\cL^{\nu}(f\zeta)-(\cL^{\nu}f)\zeta\|_p\leq
C\Big((\|\cL^\nu\zeta\|_\infty+\|\zeta\|^{1-\gamma}_\infty\|\nabla\zeta\|^\gamma_\infty)\|f\|_p
+\|\nabla\zeta\|_\infty\|f\|^{1-\beta}_p\|\nabla f\|^\beta_p\Big).\label{EE44}
\end{align}
\end{enumerate}
\el
\begin{proof}
(i). By formula (\ref{La}), we have
\begin{align}
&\cL^{\nu}( f\zeta_z)(x)-\cL^{\nu}  f(x)\cdot\zeta_z(x)-f(x)\cdot\cL^{\nu}\zeta_z(x)\no\\
&\quad=\int_{\mR^d}[f(x+y)-f(x)][\zeta_z(x+y)-\zeta_z(x)]\nu(\dif y)\no\\
&\quad=\int_{|y|\leq 1}[f(x+y)-f(x)][\zeta_z(x+y)-\zeta_z(x)]\nu(\dif y)\no\\
&\quad\quad+\int_{|y|>1}[f(x+y)-f(x)][\zeta_z(x+y)-\zeta_z(x)]\nu(\dif y)\no\\
&\quad=:I^{(1)}_z(x)+I^{(2)}_z(x).\label{EW4}
\end{align}
For $I^{(1)}_z(x)$, by Fubini's theorem, Minkowski's inequality and Lemma \ref{Le4}, we have
\begin{align*}
\int_{\mR^d}\|I^{(1)}_z\|^p_p\dif z&\leq \int_{\mR^d}
\left\|\int_{|y|\leq 1}| f(\cdot+y)-f(\cdot)|\left(\int_0^1|\nabla\zeta_z|(\cdot+sy)\dif s\right)
|y|\nu(\dif y)\right\|^p_p\dif z\\
&\leq \|\nabla\zeta\|^p_p\int_{\mR^d}\left(\int_{|y|\leq 1}| f(x+y)-f(x)|\cdot
|y|\nu(\dif y)\right)^p\dif x\\
&\leq \|\nabla\zeta\|^p_p\left(\int_{|y|\leq 1}
\| f(\cdot+y)- f(\cdot)\|_p\cdot|y|\nu(\dif y)\right)^p\\
&\leq\|\nabla\zeta\|^p_p(2\|f\|_p)^{p(1-\beta)}\|\nabla f\|_p^{p\beta}\left(\int_{|y|\leq 1}|y|^{1+\beta}\nu^{(\alpha)}(\dif y)\right)^p.
\end{align*}
For $I^{(2)}_z(x)$, we similarly have
\begin{align*}
\int_{\mR^d}\|I^{(2)}_z\|^p_p\dif z&\leq 4^p(\nu^{(\alpha)}(B^c_1))^p\|\zeta\|^p_p\|f\|^p_p.
\end{align*}
Moreover, by (\ref{EU1}) and Lemma \ref{Le6}, we also have
\begin{align*}
\int_{\mR^d}\|f\cL^{\nu}  \zeta_z\|_p^p\dif z=\|\cL^{\nu}\zeta\|_p^p\|f\|^p_p\leq C\|\zeta\|_{2,p}^p\|f\|^p_p.
\end{align*}
Combining the above calculations, we obtain (\ref{EE3}).

(ii). We have
$$
\|I^{(1)}_0\|_p\leq \|f\|_{\cH^\beta}\|\nabla\zeta\|_p\int_{|y|\leq 1}|y|^{1+\beta}\nu(\dif y)\leq
\|f\|_{\cH^\beta}\|\nabla\zeta\|_p\int_{|y|\leq 1}|y|^{1+\beta}\nu^{(\alpha)}(\dif y),
$$
and by Lemma \ref{Le4},
$$
\|I^{(2)}_0\|_p\leq \|f\|_\infty(2\|\zeta\|_p)^{1-\gamma}\|\nabla\zeta\|^\gamma_p\int_{|y|>1}|y|^\gamma\nu(\dif y)\leq
\|f\|_\infty(2\|\zeta\|_p)^{1-\gamma}\|\nabla\zeta\|^\gamma_p\int_{|y|>1}|y|^\gamma\nu^{(\alpha)}(\dif y).
$$
Estimate (\ref{EE4}) follows by (\ref{EW4}) and $\|f\cL^{\nu}  \zeta\|_p\leq \|f\|_\infty\|\cL^\nu\zeta\|_p$.
As for (\ref{EE44}), it is similar.
\end{proof}

We introduce the following notion about the non-degeneracy of $\nu^{(\alpha)}$.
\bd\label{Def}
Let $\nu^{(\alpha)}$ be a L\'evy measure with the form (\ref{Eq4}). We say that $\nu^{(\alpha)}$ is nondegenerate if
the spherical part $\Sigma$ of $\nu^{(\alpha)}$ satisfies
\begin{align}
\int_{\mS^{d-1}}|\theta_0\cdot\theta|^\alpha\Sigma(\dif\theta)\not=0,\ \ \forall\theta_0\in\mS^{d-1}.\label{Spe}
\end{align}
By the compactness of $\mS^{d-1}$ and (\ref{Es2}), the above condition
is equivalent that for some constant $\kappa_1>0$,
\begin{align}
\mathrm{Re}(\psi_{\nu^{(\alpha)}}(\xi))\geq\kappa_1|\xi|^\alpha,\ \ \xi\in\mR^d.\label{Con1}
\end{align}
\ed
\br\label{Re2}
Let $L^1_t,\cdots, L^n_t$ be $n$-independent copies of L\'evy process $L^{\nu}_t$. Write
$$
{\mathbf L}_t=(L^1_t,\cdots, L^n_t).
$$
Then ${\mathbf L}_t$ is an $nd$-dimensional L\'evy process and
the characteristic function of ${\mathbf L}_1$ is given by
$\mbox{\boldmath$\psi$}(\vec{\xi})=\psi_\nu(\xi^1)+\cdots+\psi_\nu(\xi^n)$,
where $\vec{\xi}=(\xi^1,\cdots,\xi^n)\in\mR^{nd}$ with $\xi^i\in\mR^d$.
Clearly, if
$$
\mathrm{Re}(\psi_\nu(\xi))\geq\kappa_1|\xi|^\alpha,\ \ \xi\in\mR^d,
$$
then
$$
\mathrm{Re}(\mbox{\boldmath$\psi$}(\vec{\xi}))\geq \kappa_1|\vec{\xi}|^\alpha,\ \ \vec{\xi}\in\mR^{nd}.
$$
It should be noticed that the L\'evy measure $\mbox{\boldmath$\nu$}$ of ${\mathbf L}_t$ is
very singular and has the expression
$$
\mbox{\boldmath$\nu$}(\dif\vec{x})=\nu(\dif x^1)\epsilon_{0}(\dif x^2,\cdots, \dif x^n)+\cdots
+\epsilon_{0}(\dif x^1,\cdots, \dif x^{n-1})\nu(\dif x^n),
$$
where $\vec{x}=(x^1,\cdots,x^n)\in\mR^{nd}$ with $x^i\in\mR^d$,
$\epsilon_0$ denotes the Dirac measure in $\mR^{(n-1)d}$, and the generator of ${\mathbf L}_t$ is given by
\begin{align}
\mbox{\boldmath$\cL$} f(\vec{x})=\sum_{i=1}^n\int_{\mR^d}[f(x^1,\cdots, x^i+y,\cdots, x^n)-f(\vec{x})
-y^{(\alpha)}\cdot \nabla_{x^i} f(\vec{x})]\nu(\dif y).\label{Ge}
\end{align}
\er

We need the following simple result about the smoothness of the distribution density of L\'evy process
(see \cite[Lemma 3.1]{Pr} for the symmetric case).
\bp\label{Pro1}
Let $\psi_\nu$ be defined by (\ref{Fr1}) and satisfy
\begin{align}
\mathrm{Re}(\psi_\nu(\xi))\geq \kappa_1|\xi|^\alpha,\ \ \xi\in\mR^d.\label{Con}
\end{align}
Then for each $t>0$, the law of $L^{\nu}_t$ in $\mR^d$ has a smooth density $p^{\nu}_t$
with respect to the Lebesgue measure,  and  $p^{\nu}_t\in \cap_{k\in\mN}\mW^{k,1}$. In particular, by (\ref{Ep4}),
\begin{align}
\p_t p^{\nu}_t(x)=\cL^{\nu*}p^{\nu}_t(x),\ \ (t,x)\in\mR^+\times\mR^d,\label{Ep8}
\end{align}
where $\cL^{\nu*}$ is defined by (\ref{Ad}), and $p^{\nu}_t(x)$ is also called the heat kernel of $\cL^{\nu*}$.
\ep
\begin{proof}
By (\ref{Con}) and \cite[p.190, Proposition 28.1]{Sa}, $L^{\nu}_t$ has a smooth density $p^{\nu}_t$.
Let us now prove that for each $n\in\mN$,
$\nabla^n p^{\nu}_t\in L^1(\mR^d)$. By  Fourier's transform (\ref{Fr}), one sees that
$$
p^{\nu}_t(x)=\frac{1}{(2\pi)^d}\int_{\mR^d}e^{-\mathrm{i}\xi\cdot x}e^{-t\psi_\nu(\xi)}\dif \xi.
$$
Set
$$
\phi(\xi):=\int_{|y|\leq 1}(1+\mathrm{i}\xi\cdot y-e^{\mathrm{i}\xi\cdot y})\nu(\dif y).
$$
It is easy to see that $\phi$ is a smooth complex-valued function,
and by (\ref{Con}), for any $n\in\mN$ and $j_1,\cdots, j_n\in\{1,\cdots, d\}$,
$$
\xi\to\xi_{j_1}\cdots\xi_{j_n}e^{-t\phi(\xi)}\in\cS(\mR^d),
$$
where $\xi=(\xi_1,\cdots,\xi_d)$.
Since Fourier's transform $\cF$ is a bijective and continuous linear operator from $\cS(\mR^d)$ onto itself,
there is a function $f\in\cS(\mR^d)$ such that
$$
\hat f(\xi):=\cF(f)(\xi)=\xi_{j_1}\cdots\xi_{j_n}e^{-t\phi(\xi)}.
$$
On the other hand, by L\'evy-Khintchine's representation theorem (cf. \cite[Theorem 1.2.14]{Ap}),
there is a probability measure $\mu$ on $\mR^d$ such that
$$
\hat\mu(\xi):=\int_{\mR^d}e^{\mathrm{i}\xi\cdot y}\mu(\dif y)=e^{-t(\psi_\nu-\phi)(\xi)}.
$$
Thus, by the property of Fourier's transform, we have
\begin{align*}
\p_{x_{j_1}}\cdots\p_{x_{j_n}}p^{\nu}_t(x)&=\frac{(-\mathrm{i})^n}{(2\pi)^d} \int_{\mR^d}e^{-\mathrm{i}\xi\cdot x}
(\xi_{j_1}\cdots\xi_{j_n}e^{-t\phi(\xi)})e^{-t(\psi_\nu-\phi)(\xi)}\dif \xi\\
&=\frac{(-\mathrm{i})^n}{(2\pi)^d} \int_{\mR^d}e^{-\mathrm{i}\xi\cdot x}\hat f(\xi)\hat\mu(\xi)\dif \xi
=(-\mathrm{i})^n\int_{\mR^d}f(x-y)\mu(\dif y).
\end{align*}
From this, we immediately deduce that $\nabla^n p^{\nu}_t\in L^1(\mR^d)$.
\end{proof}

Using Proposition \ref{Pro1} and Lemma \ref{Le6}, we have
the following useful estimates about the heat kernel.
\bc\label{Cor1}
Let $\nu^{(\alpha)}_i, i=1,2$ be two L\'evy measures  with the form (\ref{Eq4}),
where $\nu^{(\alpha)}_1$ is nondegenerate.
Let $\nu$ be another L\'evy measure less than $\nu^{(\alpha)}_2$.
Then, there are two indexes $\delta_1,\delta_2>1$ (depending only on $\alpha$)
and constants $C_1,C_2>0$ (depending only on $d,\alpha$, $\nu^{(\alpha)}_i$ and not on $\nu$)
such that for all $t\geq 1$,
\begin{align}
\|\nabla\cL^{\nu} p^{\nu^{(\alpha)}_1}_t\|_1&\leq C_1 t^{-\delta_1},\label{Ep6}\\
\|\p_t\cL^{\nu}p^{\nu^{(\alpha)}_1}_t\|_1&\leq C_2t^{-\delta_2}.\label{Ep5}
\end{align}
\ec
\begin{proof}
First of all, by the scaling property (\ref{Sc}) and Proposition \ref{Pro1}, we have
$$
p^{\nu^{(\alpha)}_1}_t(x)=t^{-d/\alpha}p^{\nu^{(\alpha)}_1}_1(t^{-1/\alpha}x),
$$
and for each $n\in\mN$,
\begin{align}
\int_{\mR^d}|\nabla^np^{\nu^{(\alpha)}_1}_t|(x)\dif x=t^{-n/\alpha}\int_{\mR^d}
|\nabla^np^{\nu^{(\alpha)}_1}_1|(x)\dif x\leq C t^{-n/\alpha}.\label{Ep2}
\end{align}
Estimate (\ref{Ep6}) follows from Lemma \ref{Le6} by suitable choices
of $\beta$ and $\gamma$. Notice that by (\ref{Ep8}),
$$
\p_t\cL^{\nu} p^{\nu^{(\alpha)}_1}_t(x)=\cL^{\nu}\cL^{\nu^{(\alpha)}_1*} p^{\nu^{(\alpha)}_1}_t(x).
$$
Estimate (\ref{Ep5}) follows by using Lemma \ref{Le6} twice.
\end{proof}

Now we turn to recall the classical Fefferman-Stein's theorem. Fix $\alpha\in(0,2)$.
Let $\mQ^{(\alpha)}$ be the collection of all parabolic cylinders
$$
Q_\ind:=(t_0,t_0+\ind^\alpha)\times\{x\in\mR^d: |x-x_0|\leq \ind\}.
$$
For $f\in L^1_{loc}(\mR^{d+1})$, define the Hardy-Littlewood maximal function by
$$
\cM f(t,x):=\sup_{Q\in\mQ^{(\alpha)}, (t,x)\in Q}\fint_Q|f(s,y)|\dif y\dif s,
$$
and the sharp function by
$$
f^\sharp(t,x):=\sup_{Q\in\mQ^{(\alpha)}, (t,x)\in Q}\fint_Q|f(s,y)-f_Q|\dif y\dif s,
$$
where $f_Q:=\fint_Q f(s,y)\dif y\dif s=\frac{1}{|Q|}\int_Q f(s,y)\dif y\dif s$ and $|Q|$ is the Lebesgue measure of $Q$.
One says that $f\in BMO(\mR^{d+1})$ if $f^\sharp\in L^\infty(\mR^{d+1})$. Clearly, $f\in BMO(\mR^{d+1})$ if and only if
there exists a constant $C>0$ such that for any $Q\in\mQ^{(\alpha)}$, and for some $a_Q\in\mR$,
$$
\fint_Q |f(s,y)-a_Q|\dif y\dif s\leq C.
$$
The following theorem is taken from \cite[Chapter 3]{Kr3} (see also \cite[p.148 Theorem 2]{St0}).
\bt
(Fefferman-Stein's theorem) For $p\in(1,\infty)$, there exists
a constant $C=C(p,d,\alpha)$ such that for all $f\in L^p(\mR^{d+1})$,
\begin{align}
\|f\|_p\leq C\|f^\sharp\|_p.\label{EU4}
\end{align}
\et
Using this theorem, we have
\bt\label{Th2}
For $q\in(1,\infty)$, let $\sT$ be a bounded linear operator from $L^q(\mR^{d+1})$ to $L^q(\mR^{d+1})$
and also from $L^\infty(\mR^{d+1})$ to $BMO(\mR^{d+1})$. Then for any $p\in[q,\infty)$ and $f\in L^p(\mR^{d+1})$,
$$
\|\sT f\|_p\leq C\|f\|_p,
$$
where the constant $C$ depends only on $d,p,q,\alpha$ and the norms of $\|\sT\|_{L^q\to L^q}$ and $\|\sT\|_{L^\infty\to BMO}$.
\et
\begin{proof}
Since by \cite[p.13, Theorem 1]{St},
$$
\|(\sT f)^\sharp\|_q\leq 2\|\cM\sT f\|_q\leq C\|\sT f\|_q\leq C\|\sT\|_{L^q\to L^q}\|f\|_q
$$
and
$$
\|(\sT f)^\sharp\|_\infty\leq \|\sT\|_{L^\infty\to BMO}\|f\|_\infty,
$$
by the classical Marcinkiewicz's interpolation theorem (cf. \cite{St}), we have
$$
\|\sT f\|_p\stackrel{(\ref{EU4})}{\leq} C\|(\sT f)^\sharp\|_p\leq C\|f\|_p,
$$
where $p\in[q,\infty)$.
\end{proof}

\section{A maximum principle of nonlocal parabolic equation}

In this section we fix a L\'evy measure $\nu$ less than $\nu^{(\alpha)}$ for some $\alpha\in(0,2)$,
where $\nu^{(\alpha)}$ takes the form (\ref{Eq4}), and
prove basic maximum principles for nonlocal parabolic and elliptic equations for later use.
\bl
\label{Le3}
(Maximum principle) For $T>-\infty$, let $b(t,x)$ be a bounded measurable vector field on $[T,\infty)\times\mR^d$
and $u\in C([T,\infty);C^2_b(\mR^d))$. Assume that for all $(t,x)\in[T,\infty)\times\mR^d$, $u$ satisfies
\begin{align}
u(t,x)=u(T,x)+\int^t_T\cL^{\nu} u(s,x)\dif s+\int^t_T(b\cdot\nabla u)(s,x)\dif s+\int^t_T f(s,x)\dif s.\label{EY10}
\end{align}
If $f\leq 0$, then
$$
\sup_{t\geq T}\sup_{x\in\mR^d}u(t,x)\leq \sup_{x\in\mR^d}u(T,x).
$$
In particular, the above equation admits at most one solution $u\in C([T,\infty);C^2_b(\mR^d))$.
\el
\begin{proof}
Let $\chi(x)\in[0,1]$ be a nonnegative smooth function with $\chi(x)=1$ for $|x|\leq 1$, and $\chi(x)=0$ for $|x|\geq 2$.
Set for $R>0$,
$$
\chi_R(x):=\chi(R^{-1}x),
$$
and for $\delta>0$,
$$
w^\delta_R(t,x):=\chi_R(x)u(t,x)-\delta(t-T).
$$
By (\ref{EY10}), one sees that for all $(t,x)\in[T,\infty)\times\mR^d$,
$$
w^\delta_R(t,x)=w^\delta_R(T,x)+\int^t_T\cL^{\nu}w^\delta_R(s,x)\dif s+\int^t_T(b\cdot\nabla w^\delta_R)(s,x)\dif s
+\int^t_T g_R(s,x)\dif s-\delta (t-T),
$$
where
\begin{align}
g_R:=\chi_R\cL^\nu u-\cL^\nu w_R-ub\cdot \nabla\chi_R+f\chi_R.\label{EK11}
\end{align}
For fixed $S>T$ and $\delta>0$, we want to show that for large $R$,
\begin{align}
\sup_{t\in[T,S]}\sup_{x\in\mR^d}w^\delta_R(t,x)\leq \sup_{x\in\mR^d}w^\delta_R(T,x)\leq \sup_{x\in\mR^d}u(T,x).\label{EY1}
\end{align}
If this is proven, then the result follows by firstly letting $R\to\infty$ and then $\delta\to 0$.

Below, for simplicity of notation, we drop the indexes $R$ and $\delta$.
Suppose that (\ref{EY1}) does not hold, then there exists a time $t_0\in(T,S]$ and $x_0\in \mR^d$ such that
$w$ achieves its maximum at point $(t_0,x_0)$. Thus,
\begin{align}
\nabla w(t_0,x_0)=0,\label{EW10}
\end{align}
and
\begin{align}
0&\leq\varliminf_{h\downarrow 0}\frac{1}{h}(w(t_0,x_0)-w(t_0-h,x_0))\no\\
&\leq\varlimsup_{h\downarrow 0}\frac{1}{h}\int^{t_0}_{t_0-h}\cL^\nu w(s,x_0)\dif s
+\varlimsup_{h\downarrow 0}\frac{1}{h}\int^{t_0}_{t_0-h}(b\cdot\nabla w)(s,x_0)\dif s\no\\
&+\varlimsup_{h\downarrow 0}\frac{1}{h}\int^{t_0}_{t_0-h}g(s,x_0)\dif s-\delta=:I_1+I_2+I_3-\delta.\label{EY11}
\end{align}
Since for all $y\in\mR^d$,
$$
w(t_0,x_0+y)\leq w(t_0,x_0),
$$
in view of  $w\in C([T,S];C^2_b(\mR^d))$ and by (\ref{EW10}), we have
$$
I_1=\varlimsup_{h\downarrow 0}\frac{1}{h}\int^{t_0}_{t_0-h}[\cL^\nu w(s,x_0)-\cL^\nu w(t_0,x_0)]\dif s
+\cL^\nu w(t_0,x_0)\leq 0.
$$
Similarly, for $I_2$, we have
$$
I_2=\varlimsup_{h\downarrow 0}\frac{1}{h}\int^{t_0}_{t_0-h}b(s,x_0)\cdot(\nabla w(s,x_0)-\nabla w(t_0,x_0))\dif s=0.
$$
For $I_3$, recalling (\ref{EK11}) and $f\leq 0$, by (ii) of Lemma \ref{Le1} and Lemma \ref{Le6},
we have for some $\gamma\in(0,1)$,
\begin{align*}
I_3&\leq\|\chi_R\cL^\nu u-\cL^\nu (\chi_R u)\|_\infty+\frac{\|u\|_\infty\|b\|_\infty\|\nabla\chi\|_\infty}{R}\\
&\leq\frac{C(\|u\|_\infty+\|\nabla u\|_\infty)}{R^\gamma}+\frac{\|u\|_\infty\|b\|_\infty\|\nabla\chi\|_\infty}{R},
\end{align*}
where $C$ is independent of $R$.
Choosing $R$ being sufficiently large, we obtain
$$
I_1+I_2+I_3-\delta<0,
$$
a contradiction with (\ref{EY11}). Thus, we conclude the proof of (\ref{EY1}).
\end{proof}

Similarly, we also have the following maximum principle.
\bl\label{Le5}
(Maximum principle) Assume $\lambda>0$ and $b$ is a bounded measurable vector field.
Let $u\in C^2_b(\mR^{d+1})$ (resp. $u\in C^2_b(\mR^d)$) satisfy
$$
\sL^{\nu}_{b,\lambda} u:=\p_t u-\cL^{\nu} u+(b\cdot\nabla)u+\lambda u\leq 0,
\mbox{ (resp. $(\lambda-\cL^\nu)u\leq 0$).}
$$
Then $u\leq0$. In particular, $\sL^{\nu}_{b,\lambda} u=0$ (resp. $(\lambda-\cL^\nu)u=0$)
admits at most one solution in $C^2_b(\mR^{d+1})$ (resp. $C^2_b(\mR^d)$).
\el

\bc\label{Cor2}
Let $\vartheta\in\mR^d$ and $\lambda>0$.
Then for any $p>1$,
$(\p_t-\cL^{\nu}+\vartheta\cdot\nabla+\lambda)(C^\infty_0(\mR^{d+1}))$ (resp. $(\lambda-\cL^\nu)(C^\infty_0(\mR^d))$) is dense in
$L^p(\mR^{d+1})$ (resp. $L^p(\mR^d)$).
\ec
\begin{proof}
Let $g\in L^{p/(p-1)}(\mR^{d+1})$. By Hahn-Banach's theorem, it is enough to prove that
if for all $u\in C^\infty_0(\mR^{d+1})$,
$$
\int_{\mR^{d+1}}g(t,x)\cdot(\p_t-\cL^{\nu}+\vartheta\cdot\nabla+\lambda)u(t,x)\dif x\dif t=0,
$$
then $g=0$. Since for any $(s,y)\in\mR^{d+1}$, the mapping $(t,x)\mapsto u(s+t,y+x)$
belongs to $C^\infty_0(\mR^{d+1})$. Thus, we have
$$
(\p_t-\cL^{\nu}+\vartheta\cdot\nabla+\lambda)(g\star u)=0,
$$
where $g\star u$ stands for $(s,y)\mapsto\int_{\mR^{d+1}}g(t,x)u(s+t,y+x)\dif y\dif t$.
By Lemma \ref{Le5}, $g\star u=0$ for all $u\in C^\infty_0(\mR^{d+1})$, which yields that $g=0$.
\end{proof}

\section{$L^q(\mR;L^p(\mR^d))$-maximal regularity for nonlocal parabolic equation}

Let $\vartheta\in C^\infty_b(\mR;\mR^d)$ be a time dependent vector field. For $s<t$, set
$$
\Theta_{t,s}:=\int^t_s\vartheta(r)\dif r.
$$
Let $\nu$ be a L\'evy measure and satisfy (\ref{Con}). For $f\in\cS(\mR^d)$, define
\begin{align}
\cT^{\nu}_{t,s}f(x)&:=\mE f\left(x-\Theta_{t,s}+L^{\nu}_{t-s}\right)
=\cP^\nu_{t-s}f(x-\Theta_{t,s})
=\int_{\mR^d}f(y)p^{\nu}_{t-s}\left(y-x+\Theta_{t,s}\right)\dif y.\label{Es1}
\end{align}
By (\ref{Ep8}), one has
\begin{align}
\p_t\cT^{\nu}_{t,s}f(x)&=\int_{\mR^d}f(y)\p_tp^{\nu}_{t-s}
\left(y-x+\Theta_{t,s}\right)\dif y
+\int_{\mR^d}f(y)(\vartheta_t\cdot\nabla p^{\nu}_{t-s})\left(y-x+\Theta_{t,s}\right)\dif y\no\\
&=\int_{\mR^d}f(y)(\cL^{\nu*}p^{\nu}_{t-s})\left(y-x+\Theta_{t,s}\right)\dif y
-\vartheta_t\cdot\nabla \int_{\mR^d}f(y)p^{\nu}_{t-s}
\left(y-x+\Theta_{t,s}\right)\dif y\no\\
&=\cL^{\nu}\cT^{\nu}_{t,s}f(x)-\vartheta_t\cdot\nabla\cT^{\nu}_{t,s}f(x).\label{Ep9}
\end{align}
For $\lambda\geq 0$ and $f\in\cS(\mR^{d+1})$, define
$$
u(t,x):=\int^t_{-\infty}e^{-\lambda(t-s)}\cT^{\nu}_{t,s}f(s,x)\dif s,
$$
then it is easy to check by (\ref{Ep9}) that $u\in C^\infty_b(\mR^{d+1})$ and uniquely solves
\begin{align}
\p_t u-\cL^{\nu} u+\vartheta\cdot\nabla u+\lambda u=f.\label{Eq2}
\end{align}
\br
Let $\nu_1$ and $\nu_2$ be two L\'evy measures. Let $(L^{\nu_1}_t)_{t\in\mR}$
and $(L^{\nu_2}_t)_{t\in\mR}$ be two independent L\'evy processes associated with $\nu_1$ and $\nu_2$
respectively. Then it is clear that
$$
(L^{\nu_1+\nu_2}_t)_{t\in\mR}\stackrel{(d)}{=}(L^{\nu_1}_t+L^{\nu_2}_t)_{t\in\mR}.
$$
Thus, we have
\begin{align}
\cT^{\nu_1+\nu_2}_{t,s}f(x)=\cP^{\nu_1}_{t-s}\cP^{\nu_2}_{t,s}f(x-\Theta_{t,s})
=\mE\Big(\cP^{\nu_2}_{t,s}f(x+(L^{\nu_1}_t-\Theta_{t,0})-(L^{\nu_1}_s-\Theta_{s,0})\Big).\label{Re1}
\end{align}
\er

The main aim of this section is to prove the following $L^q(\mR;L^p(\mR^d))$-regularity estimate to the above $u$
when $f\in L^q(\mR;L^p(\mR^d))$.
\bt\label{Th1}
For $\alpha\in(0,2)$, let $\nu^{(\alpha)}_i, i=1,2$ be two L\'evy measures with the form (\ref{Eq4}),
where $\nu_1^{(\alpha)}$ is nondegenerate in the sense of Definition \ref{Def}.
Let $\nu_1$ and $\nu_2$ be two L\'evy measures and satisfy  that
$$
\nu_1\geq \nu^{(\alpha)}_1,\ \ \nu_2\leq\nu_2^{(\alpha)},
$$
and for all $0<r<R<+\infty$,
$$
1_{\alpha=1}\int_{r\leq |y|\leq R}y\nu_2(\dif y)=0.
$$
Let $\vartheta:\mR\to\mR^d$ be a bounded measurable function,
and $\cT^{\nu_1}_{t,s}$ be defined by (\ref{Es1}).
Then for any $p,q\in(1,\infty)$, there exists a constant $C=C(\nu^{(\alpha)}_1,\nu^{(\alpha)}_2,\alpha,p,q,d)>0$
such that for any $-\infty\leq T<S\leq\infty$, $f\in L^q((T,S);L^p(\mR^d))$ and $\lambda\geq 0$,
\begin{align}
\int^S_T\left\|\cL^{\nu_2}\int^t_Te^{-\lambda(t-s)}
\cT^{\nu_1}_{t,s}f(s,\cdot)\dif s\right\|_p^q\dif t
\leq C\int^S_T\|f(t,\cdot)\|^q_p\dif t.\label{Es44}
\end{align}
\et
\begin{proof}
By replacing $f(t,x)$ by $f(t,x)1_{(T,S)}(t)$, it is enough to prove that
\begin{align}
\int^\infty_{-\infty}\left\|\cL^{\nu_2}\int^t_{-\infty}e^{-\lambda(t-s)}
\cT^{\nu_1}_{t,s}f(s,\cdot)\dif s\right\|_p^q\dif t
\leq C\int^\infty_{-\infty}\|f(t,\cdot)\|^q_p\dif t.\label{Es4}
\end{align}
We divide the proof into seven steps.

{\bf (Step 1)}. Let $(L^{\nu_1-\nu^{(\alpha)}_1}_t)_{t\in\mR}$ be a $d$-dimensional
L\'evy process associated with the L\'evy measure $\nu_1-\nu^{(\alpha)}_1$. By (\ref{Re1}), we have
\begin{align*}
\int^t_{-\infty}e^{-\lambda(t-s)}\cT^{\nu_1}_{t,s}f(s,x)\dif s
=\int^t_{-\infty}e^{-\lambda(t-s)}\cP^{\nu_1-\nu^{(\alpha)}_1}_{t-s}\cT^{\nu^{(\alpha)}_1}_{t,s}f(s,x)\dif s
=\mE u\Big(t,x+L^{\nu_1-\nu^{(\alpha)}_1}_t-\Theta_{t,0}\Big),
\end{align*}
where
$$
u(t,x):=\int^t_{-\infty}e^{-\lambda(t-s)}\cP^{\nu^{(\alpha)}_1}_{t-s}
f\Big(s,x-L^{\nu_1-\nu^{(\alpha)}_1}_s+\Theta_{s,0}\Big)\dif s.
$$
Suppose that (\ref{Es4}) has been proven for $\nu_1=\nu_1^{(\alpha)}$. By Fubini's theorem and Minkowski's inequality,
we have for $f\in\cS(\mR^{d+1})$,
\begin{align*}
&\int^\infty_{-\infty}\left\|\cL^{\nu_2}\int^t_{-\infty}e^{-\lambda(t-s)}
\cT^{\nu_1}_{t,s}f(s,\cdot)\dif s\right\|_p^q\dif t
=\int^\infty_{-\infty}\left\|\mE\cL^{\nu_2} u\Big(t,\cdot+L^{\nu_1-\nu^{(\alpha)}_1}_t-\Theta_{t,0}\Big)\right\|_p^q\dif t\\
&\quad\leq \int^\infty_{-\infty}\mE\left\|\cL^{\nu_2} u\Big(t,\cdot+L^{\nu_1-\nu^{(\alpha)}_1}_t-\Theta_{t,0}\Big)\right\|_p^q\dif t
\stackrel{(\ref{EU1})}{=}\mE\int^\infty_{-\infty}\left\|\cL^{\nu_2} u(t,\cdot)\right\|_p^q\dif t\leq\\
&\quad\leq C\mE\int^\infty_{-\infty}\Big\|f\Big(s,\cdot-L^{\nu_1-\nu^{(\alpha)}_1}_s+\Theta_{s,0}\Big)\Big\|^q_p\dif s
=C\int^\infty_{-\infty}\|f(s,\cdot)\|^q_p\dif s.
\end{align*}
Hence, we need only to prove (\ref{Es4}) for $\nu_1=\nu^{(\alpha)}_1$ and $\vartheta_s=0$.
Below, for simplicity of notation, we write
$$
\sL:=\cL^{\nu_2},\ \ \cL:=\cL^{\nu^{(\alpha)}_1},\ \ \cP_t:=\cP^{\nu^{(\alpha)}_1}_t,\ \
\psi_1=\psi_{\nu^{(\alpha)}_1},\ \ \psi_2=\psi_{\nu_2}.
$$

{\bf (Step 2)}. Let us firstly prove (\ref{Es4}) for $p=q=2$.
For $f\in\cS(\mR^{d+1})$, let $\hat f(s,\cdot)=\cF f(s,\cdot)$.
By (\ref{Fr}), the Fourier's transform of $\cP_tf$ is clearly given by
$$
\widehat{\cP_tf}(\xi)=e^{-\psi_1(\xi)t}\hat f(\xi).
$$
By Parseval's identity and Minkowski's inequality, we have
\begin{align*}
&\int^\infty_{-\infty}\left\|\sL
\int^t_{-\infty}e^{-\lambda(t-s)}\cP_{t-s}f(s,\cdot)\dif s\right\|_2^2\dif t\\
&\quad=\int^\infty_{-\infty}\int_{\mR^d} \left|\psi_2(\xi)\int^t_{-\infty}
e^{-\lambda(t-s)-\psi_1(\xi)(t-s)}\hat f(s,\xi)\dif s\right|^2\dif\xi\dif t\\
&\quad\stackrel{(\ref{Es22})}{\leq}\kappa_0^2\int^\infty_{-\infty}\int_{\mR^d}\left(|\xi|^\alpha\int^t_{-\infty}
e^{-\mathrm{Re}(\psi_1(\xi))(t-s)}|\hat f(s,\xi)|\dif s\right)^2\dif\xi\dif t\\
&\quad\stackrel{(\ref{Con1})}{\leq}\kappa_0^2\int^\infty_{-\infty}\int_{\mR^d}\left(|\xi|^\alpha\int^t_{-\infty}
e^{-\kappa_1|\xi|^\alpha(t-s)}|\hat f(s,\xi)|\dif s\right)^2\dif\xi\dif t\\
&\quad=\kappa_0^2\int^\infty_{-\infty}\int_{\mR^d} \left(\int^{\infty}_0
|\xi|^\alpha e^{-\kappa_1|\xi|^\alpha s}|\hat f(t-s,\xi)|\dif s\right)^2\dif\xi\dif t\\
&\quad\leq \kappa_0^2\int_{\mR^d} \left(\int^\infty_0|\xi|^\alpha e^{-\kappa_1|\xi|^\alpha s}\left(\int^\infty_{-\infty}
|\hat f(t-s,\xi)|^2\dif t\right)^{1/2}\dif s\right)^2\dif\xi\\
&\quad=\frac{\kappa_0^2}{\kappa^2_1}\int_{\mR^d}\int^\infty_{-\infty} |\hat f(t,\xi)|^2\dif t\dif\xi=
\frac{\kappa_0^2}{\kappa^2_1}\int^\infty_{-\infty}\|f(t)\|^2_2\dif t.
\end{align*}
Since $\cS(\mR^{d+1})$ is dense in $L^2(\mR^{d+1})$, (\ref{Es4}) follows for $p=q=2$.

{\bf (Step 3)}.  For $f\in L^\infty(\mR^{d+1})$, define
$$
\sT f(t,x):=\left(\sL \int^t_{-\infty}e^{-\lambda(t-s)}\cP_{t-s}f(s,\cdot)\dif s\right)(x).
$$
We want to show that
\begin{align}
\sT: L^\infty(\mR^{d+1})\to BMO(\mR^{d+1})\mbox{ is a bounded linear operator.}\label{Es3}
\end{align}
More precisely,  we want to prove that there is a constant
$C>0$ independent of $\lambda$ such that
for any $f\in L^\infty(\mR^{d+1})$ with $\|f\|_\infty\leq 1$,
and any parabolic cylinder $Q=(t_0,t_0+r^\alpha)\times B_r(x_0)$,
\begin{align}
\frac{1}{|Q|}\int_{Q}|\sT f(t,x)-a_{Q}|^2\dif x\dif t \leq C,\label{ET3}
\end{align}
where $a_{Q}$ is a constant depending on $Q$.

By shifting the origin, we may assume $t_0=0,x_0=0$. On the other hand,
by the scaling properties (\ref{Sc1}) and (\ref{EU2}), if one makes the following change in (\ref{ET3}):
$$
\nu_2(B)\to r^{\alpha}\nu_2(rB)),\ \ f(t,x)\to f(\ind^\alpha t,\ind x),\ \ \lambda\to\lambda r^\alpha,
$$
then we may further assume $\ind=1$. Thus, it suffices to prove that for any
$f\in L^\infty(\mR^{d+1})$ with $\|f\|_\infty\leq 1$,
$$
\int_{Q_1}|\sT f(t,x)-a_{Q_1}|^2\dif x\dif t \leq C,
$$
where $Q_1=(0,1)\times B_1$ and
$C=C(\nu^{(\alpha)}_1,\nu^{(\alpha)}_2,\alpha,d)$  is independent of $\nu_2$ and $\lambda$.

Following Krylov \cite{Kr1}, we now split $\sT f$ into two parts:
$$
\sT f(t,x)=\sT_1f(t,x)+\sT_2f(t,x),
$$
where for $(t,x)\in(0,1)\times B_1$,
\begin{align*}
\sT_1f(t,x)&:=\sL\left(\int^t_{-1}e^{-\lambda(t-s)}\cP_{t-s}f(s,\cdot)\dif s\right)(x),\\
\sT_2f(t,x)&:=\sL\left(\int^{-1}_{-\infty}e^{-\lambda(t-s)}\cP_{t-s}f(s,\cdot)\dif s\right)(x).
\end{align*}

{\bf (Step 4)}. In this step, we treat $\sT_1f$. Let $f_\eps(t,x):=f*\rho_\eps(t,x)$ be
the mollifying approximation of $f$, where $\rho_\eps$ is the usual mollifier in $\mR^{d+1}$. Define
\begin{align*}
u_\eps(t,x)&:=\int^t_{-1}e^{-\lambda(t-s)}\cP_{t-s}f_\eps(s,x)\dif s,\\
u(t,x)&:=\int^t_{-1}e^{-\lambda(t-s)}\cP_{t-s}f(s,x)\dif s.
\end{align*}
By definition (\ref{Es1}) and $\|f\|_\infty\leq 1$, we have
\begin{align}
|u_\eps(t,x)|\leq 2,\ \ \forall (t,x)\in[-1,1]\times\mR^d,\label{EI2}
\end{align}
and by the dominated convergence theorem,
\begin{align}
\lim_{\eps\to 0}\int^1_0\int_{B_1}|u_\eps(t,x)-u(t,x)|^2\dif x\dif t=0.\label{EI1}
\end{align}
On the other hand, by Lemma \ref{Le4}, for any $\beta\in[0,\alpha\wedge 1)$, we have
for all $t\in[-1,1], x,x'\in\mR^d$,
\begin{align}
|u_\eps(t,x)-u_\eps(t,x')|&\leq\int^t_{-1}\int_{\mR^d}|p_{t-s}(y-x)-p_{t-s}(y-x')|\dif y\dif s\no\\
&\stackrel{(\ref{EW44})}{\leq} 2^{1-\beta}\int^t_{-1}\left(|x-x'| \int_{\mR^d}|\nabla p_{t-s}(y)|\dif y\right)^\beta\dif s\no\\
&\stackrel{(\ref{Ep2})}{\leq} C|x-x'|^\beta\int^t_{-1}(t-s)^{-\beta/\alpha}\dif s\leq C|x-x'|^\beta.\label{EI22}
\end{align}
Moreover, as in the beginning of this section, since $f_\eps\in C^\infty_b(\mR^{d+1})$,
by (\ref{Ep9}) and Lemma \ref{Le3}, one sees that $u_\eps\in C^\infty_b([-1,\infty)\times\mR^{d+1})$
uniquely solves
$$
\p_t u_\eps-\cL u_\eps+\lambda u_\eps=f_\eps,\ \ u_\eps(-1,x)=0.
$$
Let $\chi$ be a nonnegative smooth function with $\chi(x)=1$ for $|x|\leq 1$ and
$\chi(x)=0$ for $|x|\geq 2$. Multiplying the above equation by $\chi$, we obtain
\begin{align}
\p_t(u_\eps\chi)&=(\cL u_\eps)\chi-\lambda u_\eps\chi+f_\eps\chi\no
=\cL(u_\eps\chi)-\lambda (u_\eps\chi)+g^\chi_\eps,
\end{align}
where
$$
g^\chi_\eps:=\chi\cL u_\eps-\cL (u_\eps\chi)+f_\eps\chi.
$$
Since $\chi$ has compact support, we have for each $t\in[0,1]$,
$$
g^\chi_\eps(t,\cdot)\in C^\infty_b(\mR^d).
$$
Thus, by Lemma \ref{Le3} again, one has the representation
$$
(u_\eps\chi)(t,x)=\int^t_{-1}e^{-\lambda(t-s)}\cP_{t-s}g^\chi_\eps(s,x)\dif s.
$$
Moreover, by (\ref{EI2}), (\ref{EI22}) and  (ii) of Lemma \ref{Le1},
$$
\int^1_{-1}\|g^\chi_\eps(t,\cdot)\|^2_2\dif t
\leq C\left(\int^1_{-1}\|\chi\cL u_\eps(t)-\cL (u_\eps(t)\chi)\|_2^2\dif t+\|\chi\|^2_2\right)\leq C.
$$
Here and below, the constant $C$ is independent of $\eps$ and $\lambda$.

As in Step 2, by Fourier's transform again, we have
\begin{align*}
\int^1_0\int_{\mR^d}|\sL(u_\eps\chi)(t,x)|^2\dif x\dif t&\leq \kappa_0^2\int^1_0\int_{\mR^d}
\left|\int^{t+1}_0|\xi|^\alpha e^{-\kappa_1|\xi|^\alpha s}|\hat g^\chi_\eps(t-s,\xi)|\dif s\right|^2\dif \xi\dif t\\
&\leq \kappa_0^2\int_{\mR^d} \left(\int^1_0|\xi|^\alpha e^{-\kappa_1|\xi|^\alpha s}\left(\int^1_{s-1}|\hat g^\chi_\eps(t-s,\xi)|^2\dif t\right)^{1/2}\dif s\right)^2\dif \xi\\
&\leq \kappa_0^2\int_{\mR^d} \left(\int^1_0|\xi|^\alpha e^{-\kappa_1|\xi|^\alpha s}\left(\int^1_{-1}|\hat g^\chi_\eps(t,\xi)|^2\dif t\right)^{1/2}\dif s\right)^2\dif \xi\\
&\leq C\int_{\mR^d}\int^1_{-1}|\hat g^\chi_\eps(t,\xi)|^2\dif t\dif \xi
=C\int^1_{-1}\|g^\chi_\eps(t,\cdot)\|^2_2\dif t\leq C.
\end{align*}
Thus, by (\ref{EI2}), (\ref{EI1}),  (\ref{EI22}) and (ii) of Lemma \ref{Le1} again, we arrive at
\begin{align*}
\int_{Q_1} |\sT_1f(t,x)|^2\dif x\dif t&=\int_{Q_1} |\sL u(t,x)|^2\dif x\dif t
\leq\sup_{\eps\in(0,1)}\int^1_0\int_{B_1}|\sL u_\eps(t,x)|^2\dif x\dif t\\
&\leq\sup_{\eps\in(0,1)}\int^1_0\int_{\mR^d}|\sL u_\eps(t,x)\cdot\chi(x)|^2\dif x\dif t\leq C.
\end{align*}

{\bf (Step 5)}. In this step, we treat $\sT_2f$ and prove that for some $a_{Q_1}\in\mR$ and
some constant $C>0$ independent of $\lambda$,
\begin{align}
\int_{Q_1}|\sT_2f(t,x)-a_{Q_1}|^2\dif x\dif t\leq C.\label{EI3}
\end{align}
Note that by (\ref{Es1}),
$$
e^{\lambda t}\sT_2f(t,x)=\int^{-1}_{-\infty}e^{\lambda s}\int_{\mR^d}f(s,y)
\sL^* p_{t-s}(y-x)\dif y\dif s=:\sT_3f(t,x).
$$
In view of $\lambda\geq 0$ and $\|f\|_\infty\leq 1$, by (\ref{Ep6}), we have
for some $\delta_1>1$ and any $(t,x)\in[0,1]\times\mR^d$,
\begin{align*}
|\nabla \sT_3f(t,x)|\leq \int^{-1}_{-\infty}\int_{\mR^d}|\nabla\sL^* p_{t-s}(y)|\dif y\dif s
\leq C\int^{-1}_{-\infty}(t-s)^{-\delta_1}\dif s\leq C,
\end{align*}
and by (\ref{Ep5}), for some $\delta_2>1$ and any $t\in[0,1]$,
\begin{align*}
|\sT_3f(t,0)-\sT_3f(0,0)|&\leq \int^{-1}_{-\infty}\int_{\mR^d}|\sL^* p_{t-s}(y)
-\sL^* p_{-s}(y)|\dif y\dif s\\
&\leq \int^{-1}_{-\infty}\int_{\mR^d}\int^t_0|\p_r\sL^* p_{r-s}(y)|\dif r\dif y\dif s\\
&\leq C\int^{-1}_{-\infty}\int^t_0(r-s)^{-\delta_2}\dif r\dif s\leq C.
\end{align*}
Hence,
$$
|\sT_3f(t,x)-\sT_3f(0,0)|\leq C,\ \ \forall(t,x)\in [0,1]\times B_1,
$$
and
$$
\int_{Q_1}|\sT_2f(t,x)-e^{-\lambda t}\sT_3f(0,0)|^2\dif x\dif t\leq C.
$$
If $\lambda=0$, we immediately have (\ref{EI3}). Now let us assume $\lambda>0$.
In this case, by Lemma \ref{Le6} and (\ref{Ep2}), we have
$$
|\sT_3f(0,0)|\leq\int^{-1}_{-\infty}e^{\lambda s}\left(\int_{\mR^d}|\sL^* p_{-s}(y)|\dif y\right)\dif s
\leq C\int^{-1}_{-\infty}e^{\lambda s}\dif s=Ce^{-\lambda}/\lambda,
$$
where $C$ is independent of $\lambda$ and $f$. So,
$$
\int_{Q_1}|(1-e^{-\lambda t})\sT_3f(0,0)|^2\dif x\dif t\leq
\frac{C}{\lambda^2}\int^1_0(1-e^{-\lambda t})^2\dif t\leq \frac{C}{3},
$$
where we have used that $1-e^{-s}\leq s$ for all $s\geq 0$. Thus, we obtain (\ref{EI3})
with $a_{Q_1}=\sT_3f(0,0)$.

{\bf (Step 6)}. Combining the above Steps 3-5, we have proven (\ref{Es3}).
By Step 2 and Theorem \ref{Th2}, we get (\ref{Es4})
for $p=q\in[2,\infty)$. As for $p=q\in(1,2)$, it follows by the following duality:
Let $g\in C^\infty_0(\mR^{d+1})$. By the integration by parts formula and the change of variables, we have
\begin{align*}
&\int^\infty_{-\infty}\int_{\mR^d} \left(\sL\int^t_{-\infty}
e^{-\lambda(t-s)}\cP_{t-s}f(s,\cdot)\dif s\right)(x)\cdot g(t,x) \dif x\dif t\\
&\qquad=\int^\infty_{-\infty}\int_{\mR^d} f(t,x)\left(\sL^*\int^t_{-\infty}e^{-\lambda(t-s)}
\cP^*_{t-s}g(s,\cdot)\dif s\right)(x) \dif x\dif t,
\end{align*}
where $\sL^*$ is the adjoint operator of $\sL$ and
$\cP^*_tg(s,x):=\mE g(s,x-L^{\nu^{(\alpha)}_1}_t)$.

{\bf (Step 7)}. For $q\not=p\in(1,\infty)$, we use a trick due to Krylov \cite{Kr2}. Clearly,
it suffices to prove that for any $T>-\infty$ and
$f\in C^\infty_0([T,\infty)\times\mR^d)$,
\begin{align}
\int^\infty_T\left\|\sL\int^t_Te^{-\lambda(t-s)}
\cP_{t-s}f(s,\cdot)\dif s\right\|_p^q\dif t
\leq C\int^\infty_T\|f(t,\cdot)\|^q_p\dif t,\label{EI7}
\end{align}
where $C$ is independent of $T$.

Set
$$
u(t,x):=\int^t_Te^{-\lambda(t-s)}\cP_{t-s}f(s,x)\dif s,\ \ w(t,x):=\sL u(t,x).
$$
By (\ref{Ep9}), one can verify that $w\in C([T,\infty);C^\infty_b(\mR^d))$ and uniquely solves
$$
\p_tw-\cL w+\lambda w=\sL f,\ \ w(T,x)=0.
$$
For $\vec{x}=(x^1,\cdots,x^n)\in\mR^{nd}$ with $x^i=(x^i_1,\cdots, x^i_d)\in\mR^d$, define
$$
U(t,\vec{x}):=w(t,x^1)\cdots w(t,x^n).
$$
Then
$$
\p_t U-\mbox{\boldmath$\cL$} U+n\lambda U=F,\ \ U(T,\vec{x})=0,
$$
where $\mbox{\boldmath$\cL$}$ is defined by (\ref{Ge}) and
$$
F(t,\vec{x})=\sum_{i=1}^n\sL_{x^i} G^i(t,\vec{x}),\ \ G^i(t,\vec{x})=f(t,x^i)\prod_{k\not= i}w(t,x^k).
$$
Here $\sL_{x^i}$ means that $\sL$ acts on the component $x^i$ of $\vec{x}$.
By the maximum principle, the unique solution $U$ can be represented by
$$
U(t,\vec{x})=\int^t_Te^{-n\lambda(t-s)}
\mbox{\boldmath$\cP$}_{t-s}F(s,\vec{x})\dif s=\sum_{i=1}^n\sL_{x^i}\int^t_Te^{-n\lambda(t-s)}
\mbox{\boldmath$\cP$}_{t-s}G^i(s,\vec{x})\dif s,
$$
where $(\mbox{\boldmath$\cP$}_t)_{t\geq 0}$ is the semigroup associated with \mbox{\boldmath$\cL$}.

Thus, by Step 6 and Minkowski's inequality, we have
\begin{align*}
\int^\infty_T\|\sL u(t)\|^{np}_p\dif t&=\int^\infty_T\|w(t)\|^{np}_p\dif t
=\int^\infty_T\int_{\mR^{nd}}|U(t,\vec{x})|^p\dif\vec{x}\dif t\\
&\leq \left(\sum_{i=1}^n\left(\int^\infty_T\int_{\mR^{nd}}\left|\sL_{x^i}\int^t_Te^{-n\lambda(t-s)}
\mbox{\boldmath$\cP$}_{t-s}G^i(s,\vec{x})\dif s\right|^p\dif\vec{x}\dif t\right)^{\frac{1}{p}}\right)^p\\
&\leq C\left(\sum_{i=1}^n\left(\int^\infty_T\int_{\mR^{nd}}|G^i(t,\vec{x})|^p\dif\vec{x}\dif t\right)^{\frac{1}{p}}\right)^p\\
&=Cn\int^\infty_T\|f(t)\|_p^p\|\sL u(t)\|^{(n-1)p}_p\dif t\\
&\leq Cn\left(\int^\infty_T\|f(t)\|_p^{np}\dif t\right)^{\frac{1}{n}}
\left(\int^\infty_T\|\sL u(t)\|^{np}_p\dif t\right)^{1-\frac{1}{n}}.
\end{align*}
From this, we get that for any $n\in\mN$ and $p>1$,
$$
\int^\infty_T\|\sL u(t)\|^{np}_p\dif t\leq (Cn)^n \int^\infty_T\|f(t)\|_p^{np}\dif t.
$$
Thus, by Marcinkiewicz's interpolation theorem (cf. \cite{St}), we get (\ref{EI7}) for any $q\geq p$.
The case $q\leq p$ follows by duality as in Step 6. The whole proof is complete.
\end{proof}

We have the following important comparison result between two different L\'evy operators.
\bt\label{Th3}
Keep the same assumptions as in Theorem \ref{Th1}. For any $p\in(1,\infty)$, there exists a constant $C>0$
such that for all $u\in\cS(\mR^d)$ and $\lambda_1,\lambda_2>0$,
\begin{align}
\|(\cL^{\nu_2}-\lambda_2)u\|_p\leq C\left(1+\frac{\lambda_2}{\lambda_1}\right)\|(\cL^{\nu_1}-\lambda_1)u\|_p.\label{EE1}
\end{align}
In particular,
\begin{align}
\|\cL^{\nu_2}u\|_p\leq C\|\cL^{\nu_1}u\|_p.\label{EE11}
\end{align}
\et
\begin{proof}
For $u\in\cS(\mR^d)$, set
$$
f:=(\cL^{\nu_1}-\lambda_1) u.
$$
By Fourier's transform, it is easy to see that
$$
u(x)=\int^\infty_0 e^{-\lambda_1 t}\cP^{\nu_1}_t f(x)\dif t.
$$
Define
$$
u_T(x):=\frac{1}{T}\int^T_0\int^t_0e^{-\lambda_1 (t-s)}\cP^{\nu_1}_{t-s} f(x)\dif s\dif t
=\int^T_0\frac{T-t}{T} e^{-\lambda_1 t}\cP^{\nu_1}_t f(x)\dif t.
$$
Then
$$
u(x)-u_T(x)=\int^\infty_T e^{-\lambda_1 t}\cP^{\nu_1}_t f(x)\dif t
+\frac{1}{T} \int^T_0te^{-\lambda_1 t}\cP^{\nu_1}_t f(x)\dif t.
$$
In view of $\|\cP^{\nu_1}_t f\|_p\leq \|f\|_p$, we have
\begin{align}
\|u-u_T\|_p\leq\|f\|_p\left(\int^\infty_T e^{-\lambda_1 t}\dif t
+\frac{1}{T}\int^\infty_0 te^{-\lambda_1 t}\dif t\right)
=\|f\|_p(\lambda_1^{-1}e^{-\lambda_1 T}+\lambda_1^{-2}T^{-1}).\label{EE2}
\end{align}
On the other hand, by (\ref{Es44}) we have
\begin{align*}
\|(\cL^{\nu_2}-\lambda_2) u_T\|^p_p&\leq\frac{1}{T}\int^T_0\left\|(\cL^{\nu_2}-\lambda_2)\int^t_0e^{-\lambda_1 (t-s)}
\cP^{\nu_1}_{t-s} f(\cdot)\dif s\right\|^p_p\dif t\\
&\leq C\|f\|^p_p+\frac{2^{p-1}}{T}\int^T_0\left(\lambda_2\int^t_0 e^{-\lambda_1(t-s)}\|f\|_p\dif s\right)^p\dif t\\
&\leq C\left(1+\frac{\lambda_2^p}{\lambda_1^p}\right)\|f\|^p_p=C\left(1+\frac{\lambda_2^p}{\lambda_1^p}\right)
\|(\cL^{\nu_1}-\lambda_1)u\|^p_p,
\end{align*}
which together with (\ref{EE2}) yields (\ref{EE1}). As for (\ref{EE11}), it follows by
firstly letting $\lambda_2\downarrow 0$ and then $\lambda_1\downarrow 0$.
\end{proof}

In the remaining part of this paper, we make the following assumption:
\begin{enumerate}[{\bf (H$^{(\alpha)}_\nu$)}]
\item Let $\nu^{(\alpha)}_i, i=1,2$ be two L\'evy measures with the form (\ref{Eq4}),
where $\nu_1^{(\alpha)}$ is nondegenerate in the sense of Definition \ref{Def}. Let $\nu$ be a L\'evy measure
satisfying (\ref{EL22}) and
$$
\nu^{(\alpha)}_1\leq\nu\leq\nu^{(\alpha)}_2.
$$
\end{enumerate}
Let $\sD^p(\cL^\nu)$ be the domain of $\cL^\nu$ in $L^p$-space, i.e.,
$$
\sD^p(\cL^\nu):=\{u\in L^p(\mR^d): \|\cL^\nu u\|_p<+\infty\}.
$$
For $\alpha\geq 0$ and $p\geq 1$, the Bessel potential space $\mH^{\alpha,p}$ is defined as the completion of
$\cS(\mR^d)$ with respect to the norm:
$$
\|f\|^\sim_{\alpha,p}:=\|(I-\Delta)^{\frac{\alpha}{2}}u\|_p\simeq \|u\|_p+\|(-\Delta)^{\frac{\alpha}{2}}u\|_p.
$$
Notice that for $k\in\mN$ and $p>1$, $\mH^{k,p}=\mW^{k,p}$ (see \cite[p135, Theorem 3]{St}).
\bc
Assume {\bf (H$^{(\alpha)}_\nu$)} with $\alpha\in(0,2)$. For any $p>1$, $f\in L^p(\mR^d)$ and
$\lambda>0$, the equation $(\cL^\nu-\lambda) u=f$ admits a unique strong solution
$u\in\mH^{\alpha,p}$. In particular, for any $p>1$, $\sD^p(\cL^\nu)=\mH^{\alpha,p}$ and
\begin{align}
\|\cL^\nu u\|_p\simeq\|(-\Delta)^{\frac{\alpha}{2}}u\|_p,\label{EW88}
\end{align}
and if $\alpha=1$, then
\begin{align}
\|\cL^\nu u\|_p\simeq\|\nabla u\|_p.\label{EW8}
\end{align}
\ec
\begin{proof}
Let $\nu^{(\alpha)}_0$ be the L\'evy measure associated with $(-\Delta)^{\frac{\alpha}{2}}$ (see (\ref{EW5})).
In Theorem \ref{Th3}, let us take $\nu_1=\nu^{(\alpha)}_0,\nu_2=\nu$ and
$\nu_1=\nu, \nu_2=\nu^{(\alpha)}_0$ respectively, then there exist $C_1,C_2>0$ such that for any
$u\in\cS(\mR^d)$ and $\lambda_1,\lambda_2>0$,
\begin{align}
\|((-\Delta)^{\frac{\alpha}{2}}+\lambda_2)u\|_p&\leq C_1\left(1+\frac{\lambda_2}{\lambda_1}\right)\|(\cL^\nu-\lambda_1)u\|_p,\label{EW6}\\
\|(\cL^\nu-\lambda_1)u\|_p&\leq C_2\left(1+\frac{\lambda_1}{\lambda_2}\right)\|((-\Delta)^{\frac{\alpha}{2}}+\lambda_2)u\|_p.\label{EW7}
\end{align}
For $\lambda>0$ and $f\in L^p(\mR^d)$, by Corollary \ref{Cor2}, there exists a sequence $u_n\in C^\infty_0(\mR^d)$
such that
$$
(\cL^\nu-\lambda) u_n\stackrel{L^p}{\to} f.
$$
By (\ref{EW6}), $u_n$ is a Cauchy sequence in $\mH^{\alpha,p}$. Let $u\in\mH^{\alpha,p}$ be the limit point.
By (\ref{EW7}), one finds that $(\cL^\nu-\lambda)u=f$. As for (\ref{EW88}), it follows by (\ref{EE11}),
and (\ref{EW8}) follows by the boundedness of Riesz transform in $L^p$-space (cf. \cite[Chapter III]{St}).
\end{proof}
\bc\label{Cor3}
Assume {\bf (H$^{(\alpha)}_\nu$)} with $\alpha\in(0,2)$. Then for any $p>1$,
$(\cP^\nu_t)_{t\geq 0}$ forms an analytic semigroup in $L^p$-space.
\ec
\begin{proof}
By \cite[Theorem 5.2]{Pa}, it suffices to prove that
$$
\|\cL^\nu\cP^\nu_t f\|_p\leq Ct^{-1}\|f\|_p,\ \ t>0,\ \  f\in L^p(\mR^d).
$$
By (\ref{Re1}), we have for any $f\in \cS(\mR^d)$,
$$
\cP^\nu_t f=\cP^{\nu^{(\alpha)}_1}_t\cP^{\nu-\nu^{(\alpha)}_1}_tf.
$$
Thus, by (\ref{Ep2}), we have
$$
\|\Delta\cP^\nu_t f\|_p\leq C t^{-\frac{2}{\alpha}}\|\cP^{\nu-\nu^{(\alpha)}_1}_tf\|_p\leq Ct^{-\frac{2}{\alpha}}\|f\|_p.
$$
Since $\cS(\mR^d)$ is dense in $L^p(\mR^d)$, we further have for any $f\in L^p(\mR^d)$,
$$
\|\Delta\cP^\nu_t f\|_p\leq Ct^{-\frac{2}{\alpha}}\|f\|_p.
$$
Now, by (\ref{EW8}) and the Gargliado-Nirenberge's inequality (cf. \cite[p.168]{Be-Lo}), we have
$$
\|\cL^\nu\cP^\nu_t f\|_p\leq C\|(-\Delta)^{\frac{\alpha}{2}}\cP^\nu_tf\|_p\leq C\|\cP^\nu_t f\|^{1-\frac{\alpha}{2}}_p
\|\Delta\cP^\nu_t f\|_p^{\frac{\alpha}{2}}\leq Ct^{-1}\|f\|_p,
$$
where $C$ is independent of $t$ and $f$.
\end{proof}

\section{Critical nonlocal parabolic equation with various coefficients}

In this section we assume {\bf (H$^{(1)}_\nu$)} with critical index $\alpha=1$. For simplicity of notation, we write
$$
\cL=\cL^\nu.
$$
Consider the following Cauchy problem of the first order critical parabolic system:
\begin{align}
\p_t u=\cL u+b\cdot\nabla u+f,\ \ u(0)=\varphi, \label{Eq3}
\end{align}
where $u=(u^1,\cdots,u^m)$, $f:\mR^+\times\mR^d\to\mR^m$, $\varphi:\mR^d\to\mR^m$ are measurable functions,
and $b:\mR^+\times\mR^d\to\mR^d$ is a bounded measurable vector field and satisfies
\begin{align}
|b(t,x)-b(t,y)|\leq\omega_b(|x-y|),\label{Es7}
\end{align}
where $\omega_b:\mR^+\to\mR^+$ is an increasing function with $\lim_{s\downarrow 0}\omega_b(s)=0$.

For obtaining the optimal regularity about the initial value, we need the following real interpolation space:
for $p>1$ and $\beta\in(0,1)$, let $\mW^{\beta,p}$ be the real interpolation space
(called Sobolev-Slobodeckij space) between $L^p$ and $\mW^{1,p}$.
By \cite[p.190,(15)]{Tr}, an equivalent norm in $\mW^{\beta,p}$ is given by
\begin{align}
\|f\|_{\beta,p}:=\|f\|_p+\left(\int_{\mR^d}\int_{\mR^d}\frac{|f(x)-f(y)|^p}
{|x-y|^{d+\beta p}}\dif x\dif y\right)^{1/p}.\label{LL1}
\end{align}
We remark that for $p\geq 2$, $\mH^{\beta,p}\subset\mW^{\beta,p}$, and for $p\leq 2$, $\mW^{\beta,p}\subset\mH^{\beta,p}$
(cf. \cite[p.155, Theorem 5 (A) and (C)]{St}). Moreover, by Sobolev's embedding theorem (see \cite[p.203, (5)]{Tr}),
if $\beta p>d$ and $\beta-\frac{d}{p}$ is not an integer, then
\begin{align}
\mW^{\beta,p}\hookrightarrow \cH^{\beta-\frac{d}{p}},\label{Sob}
\end{align}
where for $\gamma>0$, $\cH^\gamma$ is the usual H\"older space.

Let us first prove the following important apriori estimate by using the classical method of freezing coefficients (cf. \cite{Kr3}).
\bl\label{Le7}
For given $p\in(1,\infty)$, let $f\in L^p_{loc}(\mR^+;L^p(\mR^d;\mR^m))$ and
$$
u\in C(\mR^+_0; \mW^{1-\frac{1}{p},p}(\mR^d;\mR^m))\cap L^p_{loc}(\mR^+_0;\mW^{1,p}(\mR^d;\mR^m)).
$$
Assume that {\bf (H$^{(1)}_\nu$)} and (\ref{Es7}) hold,  and $u$ satisfies
\begin{align}
\p_tu(t,x)=\cL u(t,x)+b(t,x)\cdot\nabla u(t,x)+f(t,x),\ \ a.e.\ (t,x)\in\mR^+\times\mR^d.\label{EE5}
\end{align}
Then for any $T>0$,
\begin{align}
\sup_{t\in[0,T]}\|u(t)\|_{1-\frac{1}{p},p}^p+\int^T_0\|\nabla u(t)\|_p^p\dif t
\leq C(1+T^p) e^{CT^{p-1}}\left(\|u(0)\|^p_{1-\frac{1}{p},p}+\int^T_0\|f(t)\|_p^p\dif t\right),
\label{Es8}
\end{align}
where the constant $C$ depends only on $p,d$, $\|b\|_\infty$,
the modulus function $\omega_b$ and the L\'evy measures $\nu^{(1)}_i, i=1,2$.
Moreover, $u$ also satisfies the following integral equation:
\begin{align}
u(t,x)=\cP_tu(0,x)+\int^t_0\cP_{t-s}(b(s)\cdot\nabla u(s))(x)\dif s+\int^t_0\cP_{t-s}f(s,x)\dif s,\label{EQ1}
\end{align}
where $\cP_t$ is the heat semigroup associated with $\cL$.
\el
\begin{proof}
Let $(\rho_\eps)_{\eps\in(0,1)}$ be a family of mollifiers in $\mR^d$. Define
$$
u_\eps(t):=u(t)*\rho_\eps,\ \ b_\eps(t):=b(t)*\rho_\eps,\ \ f_\eps(t):=f(t)*\rho_\eps.
$$
Taking convolutions for both sides of (\ref{EE5}), we obtain
\begin{align}
\p_tu_\eps(t,x)=\cL u_\eps(t,x)+b_\eps(t,x)\cdot\nabla u_\eps(t,x)+F_\eps(t,x),\label{Eq33}
\end{align}
where
$$
F_\eps(t,x):=[(b(t)\cdot\nabla u(t))*\rho_\eps](x)-b_\eps(t,x)\cdot\nabla u_\eps(t,x)+f_\eps(t,x).
$$
Moreover, by Duhamel's formula, one sees that
\begin{align}
u_\eps(t,x)=\cP_t u_\eps(0,x)+\int^t_0\cP_{t-s}(b_\eps(s)\cdot\nabla u_\eps(s))(x)\dif s+\int^t_0\cP_{t-s}F_\eps(s,x)\dif s.\label{EQ11}
\end{align}
By the assumptions, it is easy to see that for all $\eps\in(0,1)$,
$$
|b_\eps(t,x)-b_\eps(t,y)|\leq\omega_b(|x-y|),\ \ |b_\eps(t,x)-b(t,x)|\leq\omega_b(\eps),
$$
and
$$
\lim_{\eps\to0}\int^T_0\|F_\eps(t)-f(t)\|^p_p\dif t=0.
$$
Taking limits for both sides of (\ref{EQ11}), one finds that (\ref{EQ1}) holds.
Below, we use the method of freezing the coefficients to prove
\begin{align}
\sup_{t\in[0,T]}\|u_\eps(t)\|_p^p+\int^T_0\|\nabla u_\eps(t)\|_p^p\dif t\leq
C(1+T^p) e^{CT^{p-1}}\left(\|u_\eps(0)\|^p_{1-\frac{1}{p},p}+C\int^T_0\|F_\eps(t)\|_p^p\dif t\right),\label{Ep7}
\end{align}
where the constant $C$ is independent of $\eps$ and $T$.

For simplicity of notation, we drop the subscript $\eps$ below.
Fix $\delta>0$ being small enough, whose value will be determined below.
Let $\zeta$ be a smooth function with support in $B_\delta$ and $\|\zeta\|_p=1$. For $z\in\mR^d$, set
$$
\zeta_z(x):=\zeta(x-z).
$$
Multiplying both sides of (\ref{Eq33}) by $\zeta_z$, we obtain
$$
\p_t(u\zeta_z)=(\cL u)\zeta_z+(b\cdot\nabla u)\zeta_z+F\zeta_z
=\cL(u\zeta_z)+\vartheta^b_z\cdot\nabla(u\zeta_z)+g^\zeta_z,
$$
where $\vartheta^b_z(t):=b(t,z)$ and
$$
g^\zeta_z:=(b-\vartheta^b_z)\cdot\nabla (u\zeta_z)
-ub\cdot\nabla\zeta_z+(\cL u)\zeta_z-\cL(u\zeta_z)+F\zeta_z.
$$
By Lemma \ref{Le3}, $u\zeta_z$ can be uniquely written as
$$
u\zeta_z(t,x)=\cT^{\vartheta^b_z}_{t,0}(u(0)\zeta_z)(x)+\int^t_0\cT^{\vartheta^b_z}_{t,s} g^\zeta_z(s,x)\dif s,
$$
where $\cT^{\vartheta^b_z}_{t,s}$ is defined by (\ref{Es1}) through $\vartheta^b_z$.
Thus, we have
\begin{align*}
\int^T_0\|\nabla(u\zeta_z)(t,\cdot)\|^p_p\dif t
&\leq 2^{p-1}\int^T_0\|\nabla\cT^{\vartheta^b_z}_{t,0}(u(0)\zeta_z)\|_p^p\dif t
+2^{p-1}\int^T_0\left\|\nabla\int^t_0\cT^{\vartheta^b_z}_{t,s} g^\zeta_z(s,\cdot)\dif s\right\|^p_p\dif t\\
&=:I_1(T,z)+I_2(T,z).
\end{align*}
For $I_1(T,z)$, by Corollary \ref{Cor3} and \cite[p.96 Theorem 1.14.5]{Tr}, we have
\begin{align}
\int^T_0\|\nabla\cT^{\vartheta^b_z}_{t,0}(u(0)\zeta_z)\|^p_p\dif t&\stackrel{(\ref{Es1})}{=}
\int^T_0\left\|\nabla\cP_t(u(0)\zeta_z)\left(\cdot-\int^t_0\vartheta^b_z(s)\dif s\right)\right\|^p_p\dif t
=\int^T_0\left\|\nabla\cP_t(u(0)\zeta_z)\right\|^p_p\dif t\no\\
&\stackrel{(\ref{EW8})}{\leq}C\int^T_0\left\|\cL\cP_t(u(0)\zeta_z)\right\|^p_p\dif t
\leq C\|u(0)\zeta_z\|^p_{1-\frac{1}{p},p}.\label{ET2}
\end{align}
Here and below, $C$ is independent of $T$. Thus, by definition (\ref{LL1}), it is easy to see that
$$
\int_{\mR^d}I_1(T,z)\dif z\leq C\int_{\mR^d}\|u(0)\zeta_z\|^p_{1-\frac{1}{p},p}\dif z
\leq C\Big(\|u(0)\|^p_{1-\frac{1}{p},p}\|\zeta\|_p^p+\|u(0)\|_p^p\|\zeta\|^p_{1-\frac{1}{p},p}\Big).
$$
For $I_2(T,z)$, by (\ref{EW8}) and Theorem \ref{Th1}, we have
\begin{align*}
I_2(T,z)&\leq C\int^T_0\|g^\zeta_z(s,\cdot)\|^p_p\dif s\leq
C\int^T_0\|((b-\vartheta^b_z)\cdot\nabla(u\zeta_z))(s,\cdot)\|^p_p\dif s\\
&\quad+C\int^T_0\|(ub\cdot\nabla\zeta_z)(s,\cdot)\|^p_p\dif s+C\int^T_0\|F\zeta_z(s,\cdot)\|^p_p\dif s\\
&\quad+C\int^T_0\|((\cL u)\zeta_z-\cL(u\zeta_z))(s,\cdot)\|^p_p\dif s\\
&=:I_{21}(T,z)+I_{22}(T,z)+I_{23}(T,z)+I_{24}(T,z).
\end{align*}
For $I_{21}(T,z)$, by (\ref{Es7}) and $\|\zeta\|_p=1$, we have
\begin{align*}
\int_{\mR^d}I_{21}(T,z)\dif z&\stackrel{(\ref{Es7})}{\leq}
C\omega^p_b(\delta)\int^T_0\int_{\mR^d}\|\nabla(u\zeta_z)(s,\cdot)\|^p_p\dif z\dif s\\
&\leq C\omega^p_b(\delta)\int^T_0\|\nabla u(s)\|^p_p\dif s+C\omega^p_b(\delta)\|\nabla\zeta\|^p_p\int^T_0\|u(s)\|^p_p\dif s.
\end{align*}
For $I_{24}(T,z)$, by (i) of Lemma \ref{Le1}, we have
$$
\int_{\mR^d}I_{24}(T,z)\dif z\leq C\int^T_0\|u(s)\|^p_p\dif s+C\int^T_0\|u(s)\|^{p/2}_p\|\nabla u(s)\|^{p/2}_p\dif s.
$$
Moreover, it is easy to see that
\begin{align*}
\int_{\mR^d}I_{22}(T,z)\dif z&\leq C\|b\|^p_\infty\|\nabla\zeta\|_p^p\int^T_0\|u(s)\|^p_p\dif s,\\
\int_{\mR^d}I_{23}(T,z)\dif z&\leq C\int^T_0\|F(s)\|^p_p\dif s.
\end{align*}
Combining the above calculations, we get
\begin{align*}
&\int^T_0\|\nabla u(s)\|^p_p\dif s=\int^T_0\int_{\mR^d}\|\nabla u(s)\cdot\zeta_z\|^p_p\dif z\dif s\\
&\quad\leq 2^{p-1}\int^T_0\int_{\mR^d}\|\nabla(u\zeta_z)(s)\|^p_p\dif z\dif s
+2^{p-1}\|\nabla\zeta\|_p^p\int^T_0\|u(s)\|^p_p\dif s\\
&\quad\leq C\|u(0)\|^p_{1-\frac{1}{p},p}+C\omega^p_b(\delta)\int^T_0\|\nabla u(s)\|^p_p\dif s+C\int^T_0\|u(s)\|^p_p\dif s\\
&\quad\quad+C\int^T_0\|u(s)\|^{p/2}_p\|\nabla u(s)\|^{p/2}_p\dif s+C\int^T_0\|F(s)\|^p_p\dif s.
\end{align*}
Using Young's inequality and letting $\delta$ be small enough so that $C\omega^p_b(\delta)\leq \frac{1}{4}$,
we arrive at
\begin{align}
\int^T_0\|\nabla u(s)\|^p_p\dif s\leq C\|u(0)\|^p_{1-\frac{1}{p},p}+
C\int^T_0\|u(s)\|^p_p\dif s+C\int^T_0\|F(s)\|^p_p\dif s.\label{EL1}
\end{align}
On the other hand, by (\ref{EQ11}), it is easy to see that
$$
\|u(t)\|^p_p\leq C\|u(0)\|_p^p+C t^{p-1}\|b\|^p_\infty\int^t_0 \|\nabla u(s)\|^p_p\dif s+Ct^{p-1}\int^t_0 \|F(s)\|^p_p\dif s,
$$
which together with (\ref{EL1}) and Gronwall's inequality yields that for any $T>0$,
$$
\sup_{t\in[0,T]}\|u(t)\|^p_p+\int^T_0\|\nabla u(s)\|^p_p\dif s\leq C(1+T^p) e^{CT^{p-1}}
\left(\|u(0)\|^p_{1-\frac{1}{p},p}+\int^T_0 \|F(s)\|^p_p\dif s\right).
$$
Thus, we conclude the proof of (\ref{Ep7}), and therefore,
\begin{align}
\int^T_0\|\nabla u(s)\|^p_p\dif s\leq C(1+T^p) e^{CT^{p-1}}
\left(\|u(0)\|^p_{1-\frac{1}{p},p}+\int^T_0 \|f(s)\|^p_p\dif s\right).\label{EW9}
\end{align}

Lastly, we show (\ref{Es8}). From equation (\ref{EE5}) and using estimate (\ref{EW9}), we have
\begin{align*}
\int^T_0\|\p_tu(t)\|^p_p\dif t&\leq C\left(\int^T_0\|\cL u(t)\|^p_p\dif t+\|b\|_\infty^p\int^T_0\|\nabla u(t)\|^p_p\dif t
+\int^T_0\|f(t)\|^p_p\dif t\right)\\
&\stackrel{(\ref{EW8})}{\leq} C\left((1+\|b\|_\infty^p)\int^T_0\|\nabla u(t)\|^p_p\dif t
+\int^T_0\|f(t)\|^p_p\dif t\right)\\
&\leq C(1+T^p) e^{CT^{p-1}}\left(\|u(0)\|^p_{1-\frac{1}{p},p}+\int^T_0 \|f(s)\|^p_p\dif s\right).
\end{align*}
Noticing the following embedding relation (cf. \cite[p.180, Theorem III 4.10.2]{Am})
$$
L^p([0,T],\mW^{1,p})\cap \mW^{1,p}([0,T],L^p)\hookrightarrow
C([0,T];\mW^{1-\frac{1}{p},p}),
$$
we have
\begin{align*}
\sup_{t\in[0,T]}\|u(t)\|^p_{1-\frac{1}{p},p}&\leq C\left(\int^T_0\|\p_tu(t)\|^p_p\dif t
+\int^T_0\|u(t)\|^p_{1,p}\dif t\right)\\
&\leq C(1+T^p) e^{CT^{p-1}}\left(\|u(0)\|^p_{1-\frac{1}{p},p}+\int^T_0 \|f(s)\|^p_p\dif s\right),
\end{align*}
which together with (\ref{EW9}) yields (\ref{Es8}).
\end{proof}

Before proving the existence of strong solutions to equation (\ref{Eq3}),
we recall a well-known fact (cf. \cite{Fr}, \cite{Zh4}).
\bt\label{EY2}
(Feyman-Kac formula) Let $\nu$ be a L\'evy measure and $b\in L^\infty_{loc}(\mR^+;C^\infty_b(\mR^d;\mR^d))$,
$f\in L^1_{loc}(\mR^+;\mW^\infty(\mR^d;\mR^m))$.
For any $\varphi\in\mW^\infty(\mR^d;\mR^m)$, there exists a unique $u\in C(\mR^+_0;\mW^\infty(\mR^d;\mR^m))$
satisfying
$$
u(t,x)=\varphi(x)+\int^t_0\cL^\nu u(s,x)\dif s+\int^t_0 (b\cdot\nabla u)(s,x)\dif s+\int^t_0 f(s,x)\dif s.
$$
Moreover, $u(t,x)$ can be represented by
\begin{align}
u(t,x):=\mE\varphi(X_{-t,0}(x))+\mE\left(\int^0_{-t}f(-s,X_{-t,s}(x))\dif s\right),\ t\geq 0,\label{Rep}
\end{align}
where $\{X_{t,s}(x),t\leq s\leq 0,x\in\mR^d\}$ is defined by the following SDE:
$$
X_{t,s}(x)=x+\int^s_t b(-r,X_{t,r}(x))\dif r+\int^s_t\dif L^\nu_r,\ \ t\leq s\leq 0.
$$
\et

We are now in a position to prove
\bt
Assume {\bf (H$^{(1)}_\nu$)}  and (\ref{Es7}). Let $p\in(1,\infty)$ and
$$
\varphi\in \mW^{1-\frac{1}{p},p}(\mR^d;\mR^m),\ \ \ f\in L^p_{loc}(\mR^+_0;L^p(\mR^d;\mR^m)).
$$
Then there exists a unique $u\in C(\mR^+_0;\mW^{1-\frac{1}{p},p}(\mR^d;\mR^m))\cap
L^p_{loc}(\mR^+_0;\mW^{1,p}(\mR^d;\mR^m))$ satisfying equation (\ref{EE5}).
\et
\begin{proof}
Let $b_\eps, f_\eps$ and $\varphi_\eps$ be the mollifying approximations of $b, f$ and $\varphi$:
$$
b_\eps(t,x):=b(t)*\rho_\eps(x), \ \ f_\eps(t,x):=f(t)*\rho_\eps(x),\ \ \varphi_\eps(x):=\varphi*\rho_\eps(x).
$$
By Theorem \ref{EY2}, there exists a unique $u_\eps\in C(\mR^+_0;\mW^\infty(\mR^d;\mR^m))$
satisfying the following equation:
\begin{align}
u_\eps(t,x)=\varphi_\eps(x)+\int^t_0\cL u_\eps(s,x)\dif s+\int^t_0b_\eps(s,x)\cdot\nabla u_\eps(s,x)\dif s+
\int^t_0f_\eps(s,x)\dif s.\label{ET4}
\end{align}
First of all, by Lemma \ref{Le7}, we have the following uniform estimate: for any $T>0$,
$$
\sup_{t\in[0,T]}\|u_\eps(t)\|_{1-\frac{1}{p},p}^p+
\int^T_0\|\nabla u_\eps(t)\|_p^p\dif t\leq C\left(\|\varphi\|^p_{1-\frac{1}{p},p}+\int^T_0\|f(t)\|_p^p\dif t\right),
$$
where $C$ is independent of $\eps$.

Noticing that $w_{\eps,\eps'}:=u_\eps-u_{\eps'}$ satisfies
$$
\p_t w_{\eps,\eps'}=\cL w_{\eps,\eps'}+b_\eps\cdot\nabla w_{\eps,\eps'}+
(b_\eps-b_{\eps'})\cdot\nabla u_{\eps'}+f_\eps-f_{\eps'},\ \ w_{\eps,\eps'}(0)=\varphi_\eps-\varphi_{\eps'},
$$
by Lemma \ref{Le7} again, we also have
\begin{align*}
\sup_{t\in[0,T]}\|w_{\eps,\eps'}(t)\|^p_{1-\frac{1}{p},p}+\int^T_0\|\nabla w_{\eps,\eps'}(s)\|^p_p\dif s
&\leq C\|\varphi_\eps-\varphi_{\eps'}\|_{1-\frac{1}{p},p}^p+C\int^T_0 \|f_\eps(s)-f_{\eps'}(s)\|^p_p\dif s\\
&+C\sup_{s\in[0,T]}\|b_\eps(s)-b_{\eps'}(s)\|_\infty^p\int^T_0\|\nabla u_{\eps'}(s)\|^p_p\dif s.
\end{align*}
On the other hand, by (\ref{Es7}), it is easy to see that
$$
\sup_{s\geq 0}\|b_\eps(s)-b_{\eps'}(s)\|_\infty\leq \omega_b(\eps)+\omega_b(\eps').
$$
So, for any $T>0$,
$$
\lim_{\eps,\eps'\to 0}\left(\sup_{t\in[0,T]}\|w_{\eps,\eps'}(t)\|^p_{1-\frac{1}{p},p}
+\int^T_0\|\nabla w_{\eps,\eps'}(s)\|^p_p\dif s\right)=0,
$$
and there exists a $u\in C(\mR^+_0; \mW^{1-\frac{1}{p},p}(\mR^d;\mR^m))\cap L^p_{loc}(\mR^+_0;\mW^{1,p}(\mR^d;\mR^m))$ such that
for any $T>0$,
$$
\lim_{\eps\to 0}\left(\sup_{t\in[0,T]}\|u_\eps(t)-u(t)\|^p_{1-\frac{1}{p},p}
+\int^T_0\|\nabla u_\eps(s)-\nabla u(s)\|^p_p\dif s\right)=0.
$$
By taking limits in $L^p$-space for (\ref{ET4}), one finds that for all $t\geq 0$ and almost all $x\in\mR^d$,
$$
u(t,x)=\varphi(x)+\int^t_0\cL u(s,x)\dif s+\int^t_0 b(s,x)\cdot \nabla u(s,x)\dif s+\int^t_0 f(s,x)\dif s.
$$
The existence follows. As for the uniqueness, it follows from Lemma \ref{Le7}.
\end{proof}

Now we present an application by proving a Krylov's estimate for critical diffusion process:
\begin{align}
X_t=X_0+\int^t_0b(s,X_s)\dif s+L_t.\label{SDE}
\end{align}
\bt
Assume {\bf (H$^{(1)}_\nu$)} and (\ref{Es7}). Then there exists a solution to SDE (\ref{SDE}) such that
for fixed $T_0>0$ and any $p>d+1$, stopping time $\tau$, $0\leq T\leq S\leq T_0$ and $f\in L^p([T,S]\times\mR^d)$,
\begin{align}
\mE\left(\int^{S\wedge\tau}_{T\wedge\tau} f(s, X_s)\dif s\Big|\sF_{T\wedge\tau}\right)
\leq C\|f\|_{L^p([T,S]\times\mR^d)},\label{Kry}
\end{align}
where $C$ is independent of $f$ and $\tau$.
Here, a solution to equation (\ref{SDE}) means that there exists a probability space $(\Omega,\sF,P)$
and two c\`adl\`ag stochastic processes $X_t$ and $L_t$ defined on it such that (\ref{SDE}) is satisfied, and
$L_t$ is a L\'evy process with respect to the completed filtration
$\sF_t:=\sigma^{P}\{X_s,L_s, s\leq t\}$, and whose L\'evy measure is given by $\nu$.
\et
\begin{proof}
Let $b_\eps(t,x):=b(t)*\rho_\eps(x)$ be the mollifying approximation of $b$
and let $X^\eps_t$ solve the following SDE:
\begin{align}
X^\eps_t=X_0+\int^t_0b_\eps(s,X^\eps_s)\dif s+L_t.\label{EU5}
\end{align}
It is by now standard to prove that the laws of $\{(X^\eps_t, L_t)_{t\geq 0}, \eps\in(0,1)\}$ are tight in the space of all
c\`adl\`ag functions (for example, see \cite{Zh3}).
Thus, by Skorohod's representation theorem (cf. \cite[Theorem 3.30]{Ka}), there exist a probability space still denoted by
$(\Omega,\sF,P)$ and c\`adl\`ag stochastic processes $(X^\eps_t, L^\eps_t)_{t\geq 0}$ and
$(X_t, L_t)_{t\geq 0}$ such that $(X^\eps_t, L^\eps_t)$ almost surely converges to $(X_t, L_t)$ for each $t\geq 0$,
and
$$
X^\eps_t=X^\eps_0+\int^t_0b_\eps(s,X^\eps_s)\dif s+L^\eps_t.
$$
By taking limits for equation (\ref{EU5}), it is easy to see that $(X_t,L_t)$ is a solution of SDE (\ref{SDE}).

Fix $f\in C^\infty_0(\mR^+\times\mR^d)$ and $T_0>0$.
Let $u_\eps(t,x)\in C(\mR^+_0; C^\infty_b(\mR^d))$ solve the following PDE
$$
\p_tu_\eps-\cL u_\eps-b_\eps(T_0-\cdot,\cdot)\cdot\nabla u_\eps=-f(T_0-\cdot,\cdot),\ \ u_\eps(0)=0.
$$
Set
$$
w_\eps(t,x)=u_\eps(T_0-t,x).
$$
Then
$$
\p_tw_\eps+\cL w_\eps+b\cdot\nabla w_\eps=f,\ \ w(T_0,x)=0.
$$
Let $\tau$ be any stopping time. By Ito's formula (cf. \cite[Theorem 4.4.7]{Ap}), we have
\begin{align*}
w_\eps(t,X^\eps_t)&=w(T\wedge\tau,X^\eps_{T\wedge\tau})
+\int^t_{T\wedge\tau}(\p_sw_\eps(s)+\cL w_\eps(s)+b_\eps(s)\cdot\nabla w_\eps(s))(X^\eps_s)\dif s+\mbox{ a martingale}\\
&=w(T\wedge\tau,X^\eps_{T\wedge\tau})+\int^t_{T\wedge\tau}f(s,X^\eps_s)\dif s+\mbox{ a martingale}.
\end{align*}
Taking the conditional expectations with respect to $\sF_{T\wedge\tau}$ and by the optional theorem
(cf. \cite[Theorem 6.12]{Ka}), we obtain
$$
\mE\left(\int^{S\wedge\tau}_{T\wedge\tau} f(s, X^\eps_s)\dif s\Big|\sF_{T\wedge\tau}\right)
=\mE\Big(w(S\wedge\tau,X^\eps_{S\wedge\tau})|\sF_{T\wedge\tau}\Big)-w({T\wedge\tau},X^\eps_{T\wedge\tau}).
$$
On the other hand, since
$$
|b_\eps(t,x)-b_\eps(t,y)|\leq \omega_b(|x-y|),
$$
by  (\ref{Sob}) and (\ref{Es8}), we have
$$
\sup_{t\in[T,S]}\|u_\eps\|_\infty\leq C\sup_{t\in[T,S]}\|u_\eps(t)\|_{1-\frac{1}{p},p}\leq C\|f\|_{L^p([T,S]\times\mR^d)},
$$
where the constant $C$ is independent of $\eps$. Hence,
$$
\mE\left(\int^{S\wedge\tau}_{T\wedge\tau} f(s, X^\eps_s)\dif s\Big|\sF_{T\wedge\tau}\right)\leq C\|f\|_{L^p([T,S]\times\mR^d)}.
$$
Since $f\in C^\infty_0(\mR^+\times\mR^d)$, estimate (\ref{Kry}) now follows by taking limit $\eps\to 0$.
For general $f\in L^p([T,S]\times\mR^d)$, it follows by a standard density argument.
\end{proof}

\section{Quasi-linear first order parabolic system with critical diffusion}

In this section we study the solvability of quasi-linear first order parabolic system with critical diffusions.
Let us firstly recall and extend a result of Silvestre \cite{Si2} about the H\"older estimate of advection fractional
diffusion equations.
\bt\label{Th6}
(Silvestre \cite{Si2}) Assume that $b\in L^\infty([0,1];C^\infty_b(\mR^d;\mR^d))$ and $f\in L^\infty([0,1];C^\infty_b(\mR^d))$.
For given $a>0$, let $u\in C([0,1]; C^\infty_b(\mR^d))$ satisfy that for all $(t,x)\in[0,1]\times\mR^d$,
\begin{align}
u(t,x)=u(0,x)-a\int^t_0(-\Delta)^{\frac{1}{2}}u(s,x)\dif s+\int^t_0b(s,x)\cdot\nabla u(s,x)\dif s
+\int^t_0f(s,x)\dif s.\label{EY7}
\end{align}
Then for any $\gamma\in(0,1)$, there exist a $\beta\in(0,1)$ and $C$ depending only on $d, a,\gamma$
and $\|b\|_\infty$ such that
\begin{align}
\sup_{t\in[0,1]}\|u(t)\|_{\cH^\beta}
\leq C(\|u\|_\infty+\|f\|_\infty+\|u(0)\|_{\cH^\gamma}),\label{EY9}
\end{align}
where $\|u\|_{\cH^\beta}:=\sup_{|x-y|\leq 1}|u(x)-u(y)|/|x-y|^\beta$.
\et
\begin{proof}
By \cite[Theorem 1.1]{Si2}, there exist a $\beta_0\in(0,1)$ and $C>0$ depending only on $d, a$ and $\|b\|_\infty$ such that
\begin{align}
\|u(t)\|_{\cH^{\beta_0}}\leq Ct^{-\beta_0}(\|u\|_\infty+\|f\|_\infty), \ \ t\in(0,1].\label{ER02}
\end{align}
Recall the following probabilistic representation of $u(t,x)$ (see Theorem \ref{EY2}):
\begin{align}
u(t,x)=\mE u(0,X_{-t,0}(x))+\mE\left(\int^0_{-t}f(-s,X_{-t,s}(x))\dif s\right),\ \ t\in[0,1],\label{Prob11}
\end{align}
where $\{X_{t,s}(x),-1\leq t\leq s\leq 0,x\in\mR^d\}$ is defined by the following SDE:
\begin{align}
X_{t,s}(x)=x+\int^s_t b(-r,X_{t,r}(x))\dif r+\int^s_t\dif L_r,\ \ -1\leq t\leq s\leq 0,\label{ER01}
\end{align}
where $(L_t)_{t\leq 0}$ is the L\'evy process associated with $(-\Delta)^{\frac{1}{2}}$.

By (\ref{Prob11}) and (\ref{ER01}), we have
\begin{align}
|u(t,x)-u(0,x)|&\leq \|u(0)\|_{\cH^\gamma}\mE\|X_{-t,0}(x)-x\|^\gamma+t\|f\|_\infty\no\\
&\leq\|u(0)\|_{\cH^\gamma}(t^\gamma\|b\|_\infty+\mE\|L_{-t}\|^\gamma)+t\|f\|_\infty\no\\
&\stackrel{(\ref{Sc})}{=}\|u(0)\|_{\cH^\gamma}(t^\gamma\|b\|_\infty+t^\gamma\mE\|L_{-1}\|^\gamma)+t\|f\|_\infty\no\\
&\leq t^\gamma\Big(\|u(0)\|_{\cH^\gamma}(\|b\|_\infty+\mE\|L_{-1}\|^\gamma)+\|f\|_\infty\Big).\label{ER03}
\end{align}
For given $x,y\in\mR^d$ and $t\in(0,1]$, if $t>|x-y|^{\frac{1}{2}}$, then by (\ref{ER02}) we have
$$
|u(t,x)-u(t,y)|\leq C|x-y|^{\beta_0/2}(\|u\|_\infty+\|f\|_\infty);
$$
if $t\leq|x-y|^{\frac{1}{2}}$, then by (\ref{ER03}) we have
\begin{align*}
|u(t,x)-u(t,y)|&\leq|u(t,x)-u(0,x)|+|u(t,y)-u(0,y)|+|u(0,x)-u(0,y)|\\
&\leq 2|x-y|^{\gamma/2}\Big(\|u(0)\|_{\cH^\gamma}(\|b\|_\infty+\mE\|L_{-1}\|^\gamma)+\|f\|_\infty\Big)
+|x-y|^\gamma\|u(0)\|_{\cH^\gamma}.
\end{align*}
Estimate (\ref{EY9}) now follows by taking $\beta=\min(\gamma,\beta_0)/2$.
\end{proof}

Notice that the proof of Silvestre \cite{Si2} seems strongly depend on the scale invariance of $(-\Delta)^{\frac{1}{2}}$.
Below, we use probabilistic representation (\ref{Prob11}) again to extend Silvestre's H\"older estimate to the more general L\'evy operator
(not necessary homogeneous). Consider the following L\'evy measure
$$
\nu(\dif y)=\frac{a(y)}{|y|^{d+1}}\dif y,
$$
where $a(y)$ is a measurable function on $\mR^d$ and satisfies that
$$
c_1\leq a(y)\leq c_2,
$$
and for all $0<r<R<+\infty$,
$$
\int_{r\leq |y|\leq R}\frac{ya(y)}{|y|^{d+1}}\dif y=0.
$$
Let $\cL^\nu$ be the L\'evy operator associated to $\nu$. We have
\bc\label{Cor4}
Assume that $b\in L^\infty([0,1];C^\infty_b(\mR^d;\mR^d))$ and $f\in L^\infty([0,1];C^\infty_b(\mR^d))$.
For given $\varphi\in C^\infty_b(\mR^d)$,
let $u\in C([0,1]; C^\infty_b(\mR^d))$ satisfy that for all $(t,x)\in[0,1]\times\mR^d$,
\begin{align}
u(t,x)=\varphi(x)+\int^t_0\cL^\nu u(s,x)\dif s+\int^t_0b(s,x)\cdot\nabla u(s,x)\dif s
+\int^t_0f(s,x)\dif s.\label{EY77}
\end{align}
Then for any $\gamma\in(0,1)$, there exist a $\beta\in(0,1)$ and $C$ depending only on $d, c_1,\gamma$
and $\|b\|_\infty$ such that
\begin{align}
\sup_{t\in[0,1]}\|u(t)\|_{\cH^\beta}
\leq C(\|f\|_\infty+\|\varphi\|_\infty+\|\varphi\|_{\cH^\gamma}).\label{EY99}
\end{align}
\ec
\begin{proof}

Define
$$
\nu_0(\dif y):=c_1\dif y/|y|^{d+1},\ \ \nu_1(\dif y):=\nu(\dif y)-\nu_0(\dif y)=(a(y)-c_1)\dif y/|y|^{d+1}.
$$
Let $(L^{\nu_0}_t)_{t\leq 0}$ and $(L^{\nu_1}_t)_{t\leq 0}$ be two independent $d$-dimensional L\'evy processes
with the L\'evy measures $\nu_0$ and $\nu_1$. Then we have
$$
(L^{\nu}_t)_{t\leq 0}\stackrel{(d)}{=}(L^{\nu_0}_t+L^{\nu_1}_t)_{t\leq 0}.
$$
Recall the probabilistic representation (\ref{Prob11}) of $u(t,x)$,
where $\{X_{t,s}(x),-1\leq t\leq s\leq 0,x\in\mR^d\}$ is defined by the following SDE:
$$
X_{t,s}(x)=x+\int^s_t b(-r,X_{t,r}(x))\dif r+\int^s_t\dif L^{\nu_0}_r+\int^s_t\dif L^{\nu_1}_r,\ \ -1\leq t\leq s\leq 0.
$$
Let $\mD([-1,0])$ be the space of all c\`adl\`ag functions $\ell:[-1,0]\to\mR^d$.
Below, we fix $t_0\in[0,1]$ and a path $\ell\in \mD([-1,0])$.
Let $Y_{t,s}(x,\ell_\cdot)$ solve the following SDE:
$$
Y_{t,s}(x,\ell_\cdot)=x+\int^s_t b\Big(-r,Y_{t,r}(x,\ell_\cdot)+\ell_r-\ell_{-t_0}\Big)\dif r+\int^s_t\dif L^{\nu_0}_r,\ \ -1\leq t\leq s\leq 0.
$$
By the uniqueness of solutions to SDEs, it is easy to see that
$$
X_{-t_0,s}(x)=Y_{-t_0,s}(x,L^{\nu_1}_\cdot)+L^{\nu_1}_s-L^{\nu_1}_{-t_0},\ \ -t_0\leq s\leq 0.
$$
Substituting this into (\ref{Prob11}), we get
\begin{align}
u(t_0,x)=\mE\varphi\Big(Y_{-t_0,0}(x,L^{\nu_1}_\cdot)+L^{\nu_1}_0-L^{\nu_1}_{-t_0}\Big)
+\mE\left(\int^0_{-t_0}f\Big(-s,Y_{-t_0,s}(x,L^{\nu_1}_\cdot)+L^{\nu_1}_s-L^{\nu_1}_{-t_0}\Big)\dif s\right).\label{EM1}
\end{align}
Now let us define
\begin{align}
w(t,x,\ell_\cdot):=\mE\varphi\Big(Y_{-t,0}(x,\ell_\cdot)+\ell_0-\ell_{-t_0}\Big)
+\mE\left(\int^0_{-t}f\Big(-s,Y_{-t,s}(x,\ell_\cdot)+\ell_s-\ell_{-t_0}\Big)\dif s\right).\label{EM2}
\end{align}
Using Theorem \ref{EY2} again, one sees that
$w(t,x,\ell_\cdot)$ satisfies
\begin{align*}
w(t,x,\ell_\cdot)&=\varphi(x+\ell_0-\ell_{-t_0})+\int^t_0\cL^{\nu_0}w(s,x,\ell_\cdot)\dif s\\
&+\int^t_0b(s,x+\ell_{-s}-\ell_{-t_0})\cdot\nabla w(s,x,\ell_\cdot)\dif s\\
&+\int^t_0f(s,x+\ell_{-s}-\ell_{-t_0})\dif s,
\end{align*}
where for some $a>0$,
$\cL^{\nu_0}=-a(-\Delta)^{\frac{1}{2}}$ is the L\'evy operator associated with $\nu_0$ (see (\ref{EW5})).
Thus, by Theorem \ref{Th6}, there exist a $\beta\in(0,1)$ and $C$ depending only on $d, a,\gamma$ and $\|b\|_\infty$
such that
\begin{align}
\sup_{t\in[0,1]}\|w(t,\cdot,\ell_\cdot)\|_{\cH^\beta}
&\leq C(\|w\|_\infty+\|f\|_\infty+\|\varphi\|_{\cH^\gamma})\no\\
&\stackrel{(\ref{EM2})}{\leq} C(\|f\|_\infty+\|\varphi\|_\infty+\|\varphi\|_{\cH^\gamma}).\label{EY76}
\end{align}
On the other hand, since $(L^{\nu_0}_t)_{t\leq 0}$ and $(L^{\nu_1}_t)_{t\leq 0}$ are independent,
by (\ref{EM1}) and (\ref{EM2}), we have
$$
u(t_0,x)=\mE w(t_0,x,L^{\nu_1}_\cdot).
$$
Estimate (\ref{EY99}) now follows by (\ref{EY76}).
\end{proof}

Below, for the sake of simplicity, we write
$$
\cL=\cL^\nu.
$$
Consider the following Cauchy problem of semi-linear first order parabolic system:
\begin{align}
\p_t u=\cL u+b(u)\cdot\nabla u+f(u),\ \ u(0)=\varphi, \label{Eq03}
\end{align}
where $u(t,x)=(u^1(t,x),\cdots,u^m(t,x))$, and $\varphi(x):\mR^d\to\mR^m$,
\begin{align*}
b(t,x,u)&:[0,1]\times\mR^d\times\mR^m\to\mR^d,\\
f(t,x,u)&:[0,1]\times\mR^d\times\mR^m\to\mR^m
\end{align*}
are Borel measurable functions.

We introduce the following notion about the strong solution for equation (\ref{Eq03}).
\bd\label{Def2}
Let $p>1$ and $\varphi\in \mW^{1-\frac{1}{p},p}(\mR^d;\mR^m)$. A function
$$
u\in C([0,1];\mW^{1-\frac{1}{p},p}(\mR^d;\mR^m))\cap L^p([0,1];\mW^{1,p}(\mR^d;\mR^m))
$$
is called a strong solution of equation (\ref{Eq03}) if for all $t\in[0,1]$ and almost all $x\in\mR^d$,
$$
u(t,x)=\varphi(x)+\int^t_0\cL u(s,x)\dif s+\int^t_0 b(s,x,u(s,x))\cdot\nabla u(s,x)\dif s+\int^t_0 f(s,x,u(s,x))\dif s.
$$
\ed

We firstly prove the following uniqueness of strong solutions to equation (\ref{Eq03}).
\bl
Suppose that for any $R>0$, there are two constants $C_{f,R},C_{b,R}>0$ such that for all $t\in[0,1]$,
$x,y\in\mR^d$ and $u,u'\in\mR^m$ with $|u|, |u'|\leq R$,
$$
|f(t,x,u)-f(t,x,u')|\leq C_{f,R} |u-u'|,
$$
$$
|b(t,x,u)-b(t,y,u')|\leq\omega_{b,R}(|x-y|)+C_{b,R}|u-u'|,\
$$
where $\omega_{b,R}:\mR^+\to\mR^+$ is an increasing function with $\lim_{s\downarrow 0}\omega_{b,R}(s)=0$.
Then there exists at most one strong solution in the sense of Definition \ref{Def2} provided $p>d+1$.
\el
\begin{proof}
Let $\varphi\in \mW^{1-\frac{1}{p},p}(\mR^d;\mR^m)$ and
$$
u,\tilde u\in C([0,1];\mW^{1-\frac{1}{p},p}(\mR^d;\mR^m))\cap L^p([0,1];\mW^{1,p}(\mR^d;\mR^m))
$$
be two strong solutions of equation (\ref{Eq03}) with the same initial value $\varphi$. Let
$$
w(t,x):=u(t,x)-\tilde u(t,x).
$$
Then for all $t\in[0,1]$ and almost all $x\in\mR^d$,
$$
w(t,x)=\int^t_0\cL w(s,x)\dif s+\int^t_0 b(s,x,u(s,x))\cdot\nabla w(s,x)\dif s+\int^t_0 g(s,x)\dif s,
$$
where
$$
g(t,x):=(b(t,x,u(t,x))-b(t,x,\tilde u(t,x)))\cdot\nabla\tilde u(t,x)+f(t,x,u(t,x))-f(t,x,\tilde u(t,x)).
$$
Since $u,\tilde u\in C([0,1];\mW^{1-\frac{1}{p},p}(\mR^d;\mR^m))$, by Sobolev's embedding (\ref{Sob}), for some $C>0$,
$$
\sup_{t\in[0,1]}\|u(t)\|_\infty\leq C\sup_{t\in[0,1]}\|u(t)\|_{1-\frac{1}{p},p},\ \
\sup_{t\in[0,1]}\|\tilde u(t)\|_\infty\leq C\sup_{t\in[0,1]}\|\tilde u(t)\|_{1-\frac{1}{p},p}.
$$
Let
$$
R:=C\sup_{t\in[0,1]}\|u(t)\|_{1-\frac{1}{p},p}+C\sup_{t\in[0,1]}\|\tilde u(t)\|_{1-\frac{1}{p},p},
$$
then by the assumptions, we have for all $t\in[0,1]$ and $x,y\in\mR^d$,
\begin{align*}
|b(t,x,u(t,x))-b(t,y,u(t,y))|&\leq\omega_{b,R}(|x-y|)+C_{b,R}|u(t,x)-u(t,y)|\\
&\stackrel{(\ref{Sob})}{\leq}\omega_{b,R}(|x-y|)+C\sup_{t\in[0,1]}\|u(t)\|_{1-\frac{1}{p},p}|x-y|^{1-\frac{d+1}{p}}.
\end{align*}
Thus, by Lemma \ref{Le7} and the assumptions, for all $t\in[0,1]$, we have
\begin{align}
\|w(t)\|^p_{1-\frac{1}{p},p}
\leq C\int^t_0\|g(s)\|^p_p\dif s&\leq C\int^t_0\Big(C_{b,R}^p\|\nabla\tilde u(s)\|^p_p\|w(s)\|^p_\infty
+C_{f,R}^p\|w(s)\|^p_p\Big)\dif s\no\\
&\leq C\int^t_0\Big(\|\nabla \tilde u(s)\|^p_p+1\Big)\|w(s)\|^p_{1-\frac{1}{p},p}\dif s.\label{EW3}
\end{align}
The uniqueness follows by Gronwall's inequality.
\end{proof}

We have the following existence and uniqueness of smooth solutions for equation (\ref{Eq03}).
\bt\label{Th44}
Suppose that for all $R>0$ and $j,k=0,1,2,\cdots$, there exist $C_{b,j,k, R},C_{f,j,k, R}>0$
such that for all $(t,x)\in[0,1]\times\mR^d$ and $u\in\mR^m$ with $|u|\leq R$,
\begin{align}
|\nabla^j_x\nabla^k_u b(t,x,u)|\leq C_{b,j,k,R},\ \ \ |\nabla^j_x\nabla^k_u f(t,x,u)|\leq C_{f,j,k,R},\label{WE2}
\end{align}
and there exist $\gamma_j\in\mN$, $C_{f,j}>0$ and $h_j\in (L^1\cap L^\infty)(\mR^d)$
such that for all $(t,x,u)\in[0,1]\times\mR^d\times\mR^m$,
\begin{align}
|\nabla^j_x f(t,x,u)|\leq C_{f,j}|u|^{\gamma_j}+h_j(x),\label{WE1}
\end{align}
where $\gamma_0=1$. Then for any $\varphi\in \mW^\infty(\mR^d;\mR^m)$, there exists a unique solution
$$
u\in C([0,1];\mW^\infty(\mR^d;\mR^m))
$$
to equation (\ref{Eq03}) with initial value $\varphi$. Moreover,
\begin{align}
\sup_{t\in[0,1]}\|u(t)\|_\infty\leq e^{C_{f,0}}(\|\varphi\|_\infty+\|h_0\|_\infty),\label{EU8}
\end{align}
and for any $p>d+1$,
\begin{align}
\sup_{t\in[0,1]}\|u(t)\|_{1-\frac{1}{p},p}^p+\int^1_0\|\nabla u(t)\|_p^p\dif t\leq K_p,\label{EW22}
\end{align}
where the constant $K_p$ depends only on $p,d,\nu$ and $\|\varphi\|_{1-\frac{1}{p},p}$,
$C_{f,0}$, $\|h_0\|_\infty$, $\|h_0\|_p$, $C_{b,0,0,R}$, $C_{b,0,1,R}$ and the function
\begin{align}
\omega_{b,R}(s):=\sup_{|x-y|\leq s}\sup_{t\in[0,1]}\sup_{|u|\leq R}|b(t,x,u)-b(t,y,u)|,\ \ s>0.\label{Om}
\end{align}
\et
\begin{proof}
We construct the Picardi's approximation for equation (\ref{Eq03}) as follows.
Set $u_{0}(t,x)\equiv 0$. Since for any $u\in C([0,1];\mW^\infty(\mR^d;\mR^m))$,
by (\ref{WE2}), (\ref{WE1}) and the chain rules,
\begin{align*}
(t,x)\mapsto b(t,x,u(t,x))&\in L^\infty([0,1];C^\infty_b(\mR^d;\mR^m)),\\
(t,x)\mapsto f(t,x,u(t,x))&\in L^\infty([0,1];\mW^\infty(\mR^d;\mR^m)),
\end{align*}
by Theorem \ref{EY2}, for each $n\in\mN$, there exists a unique $u_n\in C([0,1];\mW^\infty(\mR^d;\mR^m))$ solving
the following linear equation:
\begin{align}
\p_t u_n=\cL u_n+b(u_{n-1})\cdot\nabla u_n+f(u_{n-1}),\ \ u_n(0)=\varphi.\label{EY6}
\end{align}
Set
$$
\tilde u_n(t,x):=u_n(t,x)-\int^t_0\|f(s,\cdot,u_{n-1}(s,\cdot))\|_\infty\dif s,
$$
then for each $j=1,2,\cdots,m$,
$$
\p_t\tilde u^j_n-\cL\tilde u^j_n-b(u_{n-1})\cdot\nabla\tilde u^j_n=f^j(u_{n-1})-\|f_n(u_{n-1})\|_\infty\leq 0.
$$
By Lemma \ref{Le3} and (\ref{WE1}), in view of $\gamma_0=1$, we have
\begin{align*}
\|u_n(t)\|_\infty&\leq\|\tilde u_n(t)\|_\infty+\int^t_0\|f(s,\cdot,u_{n-1}(s,\cdot))\|_\infty\dif s\\
&\leq\|\tilde u_n(0)\|_\infty+\int^t_0(C_{f,0}\|u_{n-1}(s)\|_\infty+\|h_0\|_\infty)\dif s\\
&\leq\|\varphi\|_\infty+\|h_0\|_\infty+C_{f,0}\int^t_0\|u_{n-1}(s)\|_\infty\dif s,
\end{align*}
which yields by Gronwall's inequality that
\begin{align}
\sup_{t\in[0,1]}\|u_n(t)\|_\infty\leq e^{C_{f,0}}(\|\varphi\|_\infty+\|h_0\|_\infty)=:K_0.\label{EU3}
\end{align}
We mention that this $L^\infty$-estimate can be also derived by representation formula (\ref{Rep}).

Since
$$
|b(t,x,u_{n-1}(t,x))|\leq C_{b,0,0, K_0}=:K_1,
$$
by Corollary \ref{Cor4}, there exist a $\beta\in(0,1)$ and $C$ depending only on $d,\nu,p$ and $K_1$ such that
\begin{align}
&\sup_{t\in[0,1]}\|u_n(t)\|_{\cH^\beta}
\leq C(\|f(u_{n-1})\|_\infty+\|\varphi\|_\infty+\|\varphi\|_{\cH^{1-\frac{d+1}{p}}})\no\\
&\quad\stackrel{(\ref{WE1}),(\ref{EU3}),(\ref{Sob})}{\leq} C\Big(C_{f,0}K_0+\|h_0\|_\infty
+\|\varphi\|_\infty+\|\varphi\|_{1-\frac{1}{p},p}\Big)=:K_2.\label{EY999}
\end{align}
Thus, letting $\omega_{b,K_0}$ be defined by (\ref{Om}) with $R=K_0$ and using (\ref{WE2}), (\ref{EU3}), we have
\begin{align}
|b(t,x,u_{n-1}(t,x))-b(t,y,u_{n-1}(t,y))|\leq \omega_{b,K_0}(|x-y|)+C_{b,0,1,K_0} K_2|x-y|^\beta.\label{WQ4}
\end{align}
Hence, we can use Lemma \ref{Le7} to derive that for any $p>1$,
\begin{align}
\|u_n(t)\|_{1-\frac{1}{p},p}^p+
\int^t_0\|\nabla u_n(s)\|_p^p\dif s&\leq C\left(\|\varphi\|^p_{1-\frac{1}{p},p}
+\int^t_0\|f(s,u_{n-1}(s))\|_p^p\dif s\right)\no\\
&\leq C_1\left(\|\varphi\|^p_{1-\frac{1}{p},p}+
\int^t_0\Big(C_{f,0}^p\|u_{n-1}(s)\|_p^p+\|h_0\|^p_p\Big)\dif s\right),\label{WE4}
\end{align}
where $C_1\geq 1$ depends only on $p,d,\nu$, $K_1$, $K_2$, $\omega_{b,K_0}$ and $C_{b,0,1,K_0}$.
In particular, for any $t\in[0,1]$,
$$
\|u_n(t)\|_p^p\leq C_1\left(\|\varphi\|^p_{1-\frac{1}{p},p}+\|h_0\|^p_p\right)+
C_1C_{f,0}^p\int^t_0\|u_{n-1}(s)\|_p^p\dif s,
$$
and by Gronwall's inequality,
$$
\sup_{t\in[0,1]}\|u_n(t)\|_p^p\leq C_1\left(\|\varphi\|^p_{1-\frac{1}{p},p}+\|h_0\|^p_p\right) e^{C_1C_{f,0}^p}.
$$
Substituting this into (\ref{WE4}), we obtain
\begin{align}
\sup_{t\in[0,1]}\|u_n(t)\|_{1-\frac{1}{p},p}^p+\int^1_0\|\nabla u_n(t)\|_p^p\dif t
\leq C_1\left(\|\varphi\|^p_{1-\frac{1}{p},p}
+\int^1_0\|f(s,u_{n-1}(s))\|_p^p\dif s\right)\leq K_3,\label{EW2}
\end{align}
where $K_3$ depends only on $p$, $C_1$, $\|\varphi\|_{1-\frac{1}{p},p}$, $C_{f,0}$, $\|h_0\|_p$.

Let us now estimate the higher order derivatives of $u_n$.
For given $k\in\mN$, set
$$
w^{(k)}_n(t,x):=\nabla^k u_n(t,x).
$$
By equation (\ref{EY6}) and the chain rules, one sees that
$$
\p_tw^{(k)}_n=\cL w^{(k)}_n+b(u_{n-1})\cdot\nabla w^{(k)}_n+g^{(k)}_n,
$$
where
$$
g^{(k)}_n(t,x):=\nabla^k (f(t,\cdot,u_{n-1}(t,\cdot)))(x)+\sum_{j=1}^k \frac{k!}{(k-j)!j!}\nabla^j(b(t,\cdot,u_{n-1}(t,\cdot)))(x)
\cdot\nabla^{k-j}\nabla u_n(t,x).
$$
By (\ref{WQ4}) and Lemma \ref{Le7}, for any $p>1$, we have
$$
\sup_{t\in[0,1]}\|w^{(k)}_n(t)\|_{1-\frac{1}{p},p}^p+
\int^1_0\|\nabla w^{(k)}_n(s)\|_p^p\dif s\leq C\left(\|\nabla^k\varphi\|^p_{1-\frac{1}{p},p}
+\int^1_0\|g^{(k)}_n(s)\|_p^p\dif s\right).
$$
Since $g^{(k)}_n(s)$ contains at most $k$-order derivatives of $u_n(s)$ and
the powers of lower order derivatives of $u_n(s)$,
by induction method, it is easy to see that for any $k\in\mN$ and $p>1$,
\begin{align}
\sup_{t\in[0,1]}\|w^{(k)}_n(t)\|_{1-\frac{1}{p},p}^p+
\int^1_0\|\nabla w^{(k)}_n(s)\|_p^p\dif s\leq K_{p,k},\label{WE3}
\end{align}
where $K_{p,k}$ is independent of $n$.

Define
$$
w_{n,m}(t,x):=u_n(t,x)-u_m(t,x).
$$
Then
$$
\p_t w_{n,m}=\cL w_{n,m}+b(u_{n-1})\cdot\nabla w_{n,m}
+(G_{1,n,m}+G_{2,n,m})w_{n-1,m-1},
$$
subject to $w_{n,m}(0)=0$, where
\begin{align*}
G^{ki}_{1,n,m}(t,x)&:=\sum_j\int^1_0 \p_{u^i}b^j(t,x, u_{n-1}(t,x)
+r (u_{n-1}-u_{m-1})(t,x))\dif r\cdot \p_ju^k_m(t,x),\\
G^{ki}_{2,n,m}(t,x)&:=\int^1_0 \p_{u^i}f^k(t,x, u_{n-1}(t,x)+r (u_{n-1}-u_{m-1})(t,x))\dif r.
\end{align*}
By (\ref{WQ4}) and Lemma \ref{Le7} again, we have
$$
\|w_{n,m}(t)\|^p_{1-\frac{1}{p},p}
\leq C\int^t_0\|(G_{1,n,m}(s)+G_{2,n,m}(s))w_{n-1,m-1}(s)\|^p_p\dif s.
$$
By (\ref{WE2}) and as in estimating (\ref{EW3}), we further have
\begin{align*}
\|w_{n,m}(t)\|^p_{1-\frac{1}{p},p}
&\leq C\int^t_0\Big(\|\nabla u_m(s)\|^p_p+1\Big)\|w_{n-1,m-1}(s)\|^p_{1-\frac{1}{p},p}\dif s\\
&\stackrel{(\ref{WE3})}{\leq} C(K_{p,1}+1)\int^t_0\|w_{n-1,m-1}(s)\|^p_{1-\frac{1}{p},p}\dif s.
\end{align*}
Taking super-limit for both sides and by Fatou's lemma, we obtain
$$
\varlimsup_{n,m\to\infty}\sup_{s\in[0,t]}\|w_{n,m}(s)\|^p_{1-\frac{1}{p},p}
\leq C(K_{p,1}+1)\int^t_0\varlimsup_{n,m\to\infty}\sup_{s\in[0,r]}\|w_{n-1,m-1}(s)\|^p_{1-\frac{1}{p},p}\dif r.
$$
Thus, by Gronwall's inequality, we get
$$
\varlimsup_{n,m\to\infty}\sup_{t\in[0,1]}\|w_{n,m}(t)\|^p_{1-\frac{1}{p},p}=0,
$$
which together with (\ref{WE3}) and the interpolation inequality yields that for any $k\in\mN$,
$$
\lim_{n,m\to\infty}\sup_{t\in[0,1]}\|u_n(t)-u_m(t)\|^p_{k,p}=0.
$$
Hence, there exists a $u\in C([0,1];\mW^\infty(\mR^d;\mR^m))$ such that for any $k\in\mN$,
$$
\lim_{n\to\infty}\sup_{t\in[0,1]}\|u_n(t)-u(t)\|^p_{k,p}=0.
$$
The proof is finished by taking limits for equation (\ref{EY6}).
\end{proof}

Next we show the well-posedness of equation (\ref{Eq03}) under less regularity conditions on $b,f$.
\bt\label{Th4}
Let $p>d+1$. Suppose that there exist $C_f>0$ and $h\in (L^p\cap L^\infty)(\mR^d)$ such that
for all $(t,x,u)\in [0,1]\times\mR^d\times\mR^m$,
\begin{align}
|f(t,x,u)|\leq C_f|u|+h(x);\label{WE5}
\end{align}
and for any $R>0$, there are three constants $C_{f,R},C_{b,0,R},C_{b,1,R}>0$ such that for all $t\in[0,1]$,
$x,y\in\mR^d$ and $u,u'\in\mR^m$ with $|u|, |u'|\leq R$,
\begin{align}
\left\{
\begin{aligned}
|f(t,x,u)-f(t,x,u')|&\leq C_{f,R} |u-u'|,\ \ |b(t,x,u)|\leq C_{b,0,R},\\
|b(t,x,u)-b(t,y,u')|&\leq\omega_{b,R}(|x-y|)+C_{b,1,R}|u-u'|,
\end{aligned}
\right.\label{WE6}
\end{align}
where $\omega_{b,R}:\mR^+\to\mR^+$ is an increasing function with $\lim_{s\downarrow 0}\omega_{b,R}(s)=0$.
Then for any $\varphi\in \mW^{1-\frac{1}{p},p}(\mR^d;\mR^m)$,
there exists a unique strong solution $u$ in the sense of Definition \ref{Def2}. Moreover,
\begin{align}
\sup_{t\in[0,1]}\|u(t)\|_\infty\leq e^{C_f}(\|\varphi\|_\infty+\|h\|_\infty).\label{EU88}
\end{align}
\et
\begin{proof}
We divide the proof into three steps.

{\bf (Step 1)}.
Let $\chi(x)\in[0,1]$ be a nonnegative smooth function with $\chi(x)=1$ for $|x|\leq 1$ and $\chi(x)=0$ for $|x|>2$.
Let $(\rho^x_\eps)_{\eps\in(0,1)}$ and $(\rho^u_\eps)_{\eps\in(0,1)}$ be the mollifiers in $\mR^d$ and $\mR^m$. Define
$$
b_\eps(t,x,u):=b(t,\cdot,\cdot)*(\rho^x_\eps\rho^u_\eps)(x,u),\ \ \varphi_\eps(x):=\varphi*\rho^x_\eps(x),
$$
and
$$
f_\eps(t,x,u):=f(t,\cdot,\cdot)*(\rho^x_\eps\rho^u_\eps)(x,u)\chi(\eps x).
$$
By (\ref{WE5}) and  (\ref{WE6}), one sees that (\ref{WE2}) and (\ref{WE1}) are satisfied for $b_\eps$ and $f_\eps$, and
\begin{align}
|f_\eps(t,x,u)|&\leq \Big(C_f(|u|+\eps)+h*\rho^x_\eps(x)\Big)\chi(\eps x)\no\\
&\leq C_f|u|+C_f\eps\chi(\eps x)+h*\rho^x_\eps(x),\label{WE8}
\end{align}
and for any $R>0$ and all $t\in[0,1]$, $x,y\in\mR^d$ and $u,u'\in\mR^m$ with $|u|, |u'|\leq R$,
\begin{align}
\left\{
\begin{aligned}
|f_\eps(t,x,u)-f_\eps(t,x,u')|&\leq C_{f,R+1}|u-u'|,\ \ |b_\eps(t,x,u)|\leq C_{b,R+1},\\
|b_\eps(t,x,u)-b_\eps(t,y,u')|&\leq \omega_{b,R+1}(|x-y|)+C_{b,R+1}|u-u'|.
\end{aligned}
\right.\label{EU6}
\end{align}
Moreover, by definition (\ref{LL1}),
\begin{align}
\|\varphi_\eps\|_{1-\frac{1}{p},p}\leq \|\varphi\|_{1-\frac{1}{p},p}.\label{EU66}
\end{align}
By Theorem \ref{Th44}, let $u_\eps\in C([0,1];\mW^\infty(\mR^d;\mR^m))$ solve the following equation
\begin{align}
\p_t u_\eps=\cL u_\eps+b_\eps(u_\eps)\cdot\nabla u_\eps+f_\eps(u_\eps),\ \ u_\eps(0)=\varphi_\eps.\label{EY66}
\end{align}
By (\ref{EU8}) and (\ref{WE8}), we have
\begin{align}
\sup_{t\in[0,1]}\|u_\eps(t)\|_\infty\leq e^{C_f}(\|\varphi\|_\infty+C_f\eps+\|h\|_\infty),\label{EU333}
\end{align}
and by (\ref{WE8}), (\ref{EU6}), (\ref{EU66}) and (\ref{EW22}),
\begin{align}
\sup_{\eps\in(0,1)}\left(\sup_{t\in[0,1]}\|u_\eps(t)\|_{1-\frac{1}{p},p}^p
+\int^1_0\|\nabla u_\eps(t)\|_p^p\dif t\right)\leq K,\label{EW222}
\end{align}
where we have particularly used that for $p>d+1$,
$$
\|C_f\eps\chi(\eps \cdot)+h*\rho^x_\eps\|_p\leq C_f\eps^{1-d/p}\|\chi\|_p+\|h\|_p\leq C_f\|\chi\|_p+\|h\|_p.
$$

{\bf (Step 2)}. In this step we want to show that
\begin{align}
\lim_{N\to\infty}\sup_{\eps\in(0,1)}\sup_{t\in[0,1]}\int_{|x|\geq N}|u_\eps(t,x)|^p\dif x=0.\label{WQ2}
\end{align}
Let $\zeta_N(x):=1-\chi(N^{-1}x)$. Multiplying both sides of equation (\ref{EY66}) by $\zeta_N(x)$, we have
$$
\p_t(u_\eps\zeta_N)=\cL(u_\eps\zeta_N)+b_\eps(u_\eps)\cdot\nabla(u_\eps\zeta_N)+g_{N,\eps},
$$
where
$$
g_{N,\eps}:=\zeta_N\cL u_\eps-\cL(u_\eps\zeta_N)-u_\eps b_\eps(u_\eps)\cdot\nabla\zeta_N+f_\eps(u_\eps)\zeta_N.
$$
Let
$$
R:=e^{C_f}(\|\varphi\|_\infty+C_f+\|h\|_\infty).
$$
Since
\begin{align}
|b_\eps(t,x,u_\eps(t,x))-b_\eps(t,y,u_\eps(t,y))|
&\stackrel{(\ref{EU6})}{\leq}\omega_{b,R+1}(|x-y|)+C_{b,R+1}|u_\eps(t,x)-u_\eps(t,y)|\no\\
&\stackrel{(\ref{Sob})}{\leq}\omega_{b,R+1}(|x-y|)+C\sup_{t\in[0,1]}\|u_\eps(t)\|_{1-\frac{1}{p},p}|x-y|^{1-\frac{d+1}{p}}\no\\
&\stackrel{(\ref{EW222})}{\leq}\omega_{b,R+1}(|x-y|)+CK^{\frac{1}{p}}|x-y|^{1-\frac{d+1}{p}},\label{WQ1}
\end{align}
here and below, the constant $C$ is independent of $N$ and $\eps$, by Lemma \ref{Le7}, we have
\begin{align}
\|u_\eps(t)\zeta_N\|^p_{1-\frac{1}{p},p}\leq C\|\varphi_\eps\zeta_N\|^p_{1-\frac{1}{p},p}
+C\int^t_0\|g_{N,\eps}(s)\|^p_p\dif s.\label{WQ3}
\end{align}
Clearly,
\begin{align*}
\|\varphi_\eps\zeta_N\|^p_{1-\frac{1}{p},p}\leq
C\|\varphi_\eps\zeta_N\|^p_{1,p}\leq C\|\varphi\zeta_N\|_p^p+
C\|\nabla\varphi\zeta_N\|_p^p+C\|\varphi\nabla\zeta_N\|_p^p\to 0,\ \ N\to\infty.
\end{align*}
By (\ref{EE44}) and (\ref{WE8}), we have
\begin{align*}
\|g_{N,\eps}\|_p&\leq\|\zeta_N\cL u_\eps-\cL(u_\eps\zeta_N)\|_p
+\|u_\eps b_\eps(u_\eps)\cdot\nabla\zeta_N\|_p+\|f_\eps(u_\eps)\zeta_N\|_p\\
&\leq C\Big((\|\cL\zeta_N\|_\infty+\|\zeta_N\|^{\frac{1}{2}}_\infty\|\nabla\zeta_N\|^{\frac{1}{2}}_\infty)\|u_\eps\|_p
+\|\nabla\zeta_N\|_\infty\|u_\eps\|^{\frac{1}{2}}_p\|\nabla u_\eps\|^{\frac{1}{2}}_p\Big)\\
&\quad+\|u_\eps\|_p\|b_\eps(u_\eps)\|_\infty\|\nabla\zeta_N\|_\infty+C_f\|u_\eps\zeta_N\|_p
+C_f\eps\|\chi(\eps \cdot)\zeta_N\|_p+\|(h*\rho^x_\eps)\zeta_N\|_p.
\end{align*}
Noticing that
$$
\eps^p\|\chi(\eps \cdot)\zeta_N\|^p_p=\eps^{p-d}\int_{\mR^d}|\chi(x)(1-\chi(N^{-1}\eps^{-1}x))|^p\dif x
\leq (\frac{2}{N})^{p-d}\int_{\mR^d}|\chi(x)|^p\dif x
$$
and
$$
\|(h*\rho^x_\eps)\zeta_N\|^p_p\leq \int_{B^c_{N-1}}|h(x)|^p\dif x,
$$
by Lemma \ref{Le6} and (\ref{EW222}), we have
\begin{align*}
\int^t_0\|g_{N,\eps}(s)\|^p_p\dif s\leq \frac{C}{N^{\frac{p}{2}}}+C\int^t_0\|u_\eps(s)\zeta_N\|_p^p\dif s
+\frac{C}{N^{p-d}}\int_{\mR^d}|\chi(x)|^p\dif x+C\int_{B^c_{N-1}}|h(x)|^p\dif x.
\end{align*}
Substituting this into (\ref{WQ3}) and using Gronwall's inequality, we obtain
$$
\lim_{N\to\infty}\sup_{\eps\in(0,1)}\sup_{t\in[0,1]}\|u_\eps(t)\zeta_N\|^p_p=0.
$$
This clearly implies (\ref{WQ2}).

{\bf (Step 3)}. For fixed $\eps,\eps'\in(0,1)$, let us define
$$
w_{\eps,\eps'}(t,x):=u_\eps(t,x)-u_{\eps'}(t,x).
$$
Then
$$
\p_t w_{\eps,\eps'}=\cL w_{\eps,\eps'}+b_\eps(u_\eps)\cdot\nabla w_{\eps,\eps'}
+(G_{1,\eps,\eps'}+G_{2,\eps,\eps'})w_{\eps,\eps'}+F_{1,\eps,\eps'}+F_{2,\eps,\eps'},
$$
subject to $w_{\eps,\eps'}(0)=\varphi_\eps-\varphi_{\eps'}$, where
\begin{align*}
G^{ki}_{1,\eps,\eps'}(t,x)&:=\sum_j\int^1_0 \p_{u^i}b_\eps^j(t,x, u_\eps(t,x)
+r (u_{\eps}-u_{\eps'})(t,x))\dif r\cdot \p_ju^k_{\eps'}(t,x),\\
G^{ki}_{2,\eps,\eps'}(t,x)&:=\int^1_0 \p_{u^i}f^k_\eps(t,x, u_\eps(t,x)+r (u_{\eps}-u_{\eps'})(t,x))\dif r,\\
F_{1,\eps,\eps'}(t,x)&:=(b_\eps(t,x,u_{\eps'}(t,x))-b_{\eps'}(t,x,u_{\eps'}(t,x)))\cdot\nabla u_{\eps'}(t,x),\\
F_{2,\eps,\eps'}(t,x)&:=f_\eps(t,x,u_{\eps'}(t,x))-f_{\eps'}(t,x,u_{\eps'}(t,x)).
\end{align*}
By (\ref{WQ1}) and Lemma \ref{Le7} again, we have
$$
\|w_{\eps,\eps'}(t)\|^p_{1-\frac{1}{p},p}+\int^t_0\|\nabla w_{\eps,\eps'}(s)\|^p_p\dif s
\leq h_{\eps,\eps'}+C\int^t_0\|(G_{1,\eps,\eps'}(s)+G_{2,\eps,\eps'}(s))w_{\eps,\eps'}(s)\|^p_p\dif s,
$$
where
\begin{align*}
h_{\eps,\eps'}:=C\|w_{\eps,\eps'}(0)\|_{1-\frac{1}{p},p}^p+C\int^1_0\|F_{1,\eps,\eps'}(s)+F_{2,\eps,\eps'}(s)\|^p_p\dif s.
\end{align*}
By (\ref{EU6}) and as in estimating (\ref{EW3}), we further have
$$
\|w_{\eps,\eps'}(t)\|^p_{1-\frac{1}{p},p}+\int^t_0\|\nabla w_{\eps,\eps'}(s)\|^p_p\dif s
\leq h_{\eps,\eps'}+C\int^t_0\Big(\|\nabla u_{\eps'}(s)\|^p_p+1\Big)\|w_{\eps,\eps'}(s)\|^p_{1-\frac{1}{p},p}\dif s.
$$
By Gronwall's inequality and (\ref{EW222}), one sees that
\begin{align}
\sup_{s\in[0,1]}\|w_{\eps,\eps'}(s)\|^p_{1-\frac{1}{p},p}
+\int^1_0\|\nabla w_{\eps,\eps'}(s)\|^p_p\dif s\leq Ch_{\eps,\eps'}.\label{EEQ1}
\end{align}
Now it is standard to show that
$$
\lim_{\eps,\eps'\to 0}\|w_{\eps,\eps'}(0)\|_{1-\frac{1}{p},p}^p\leq C\lim_{\eps,\eps'\to 0}\|w_{\eps,\eps'}(0)\|_{1,p}^p=0,
$$
and by (\ref{WE6}) and (\ref{EW222}),
$$
\lim_{\eps,\eps'\to 0}\int^1_0\|F_{1,\eps,\eps'}(s)\|^p_p\dif s\leq K\lim_{\eps,\eps'\to 0}
\Big(\omega_{b,R+1}(\eps)+C_{b,1,R+1}\eps+\omega_{b,R+1}(\eps')+C_{b,1,R+1}\eps'\Big)^p=0.
$$
We now look at $F_{2,\eps,\eps'}$. For any $N>0$, we write
$$
\int^1_0\!\!\!\int_{\mR^d}|F_{1,\eps,\eps'}(s,x)|^p_p\dif x\dif s
=\int^1_0\!\!\!\int_{B^c_N}|F_{1,\eps,\eps'}(s,x)|^p_p\dif x\dif s
+\int^1_0\!\!\!\int_{B_N}|F_{1,\eps,\eps'}(s,x)|^p_p\dif x\dif s=:I_1+I_2.
$$
For $I_1$, by (\ref{WE8}) we have
\begin{align*}
I_1&\leq\int^1_0\!\!\!\int_{B^c_N}\Big(2C_f|u_{\eps'}(s,x)|+\eps\chi(\eps x)
+h*\rho_\eps(x)+\eps'\chi(\eps' x)+h*\rho_{\eps'}(x)\Big)^p\dif x\dif s\\
&\leq C\sup_{s\in[0,1]}\int_{B^c_N}|u_{\eps'}(s,x)|^p\dif x
+\frac{C}{N^{p-d}}\int_{\mR^d}|\chi(x)|^p\dif x+C\int_{B^c_{N-1}}|h(x)|^p\dif x,
\end{align*}
which converges to zero uniformly in $\eps'\in(0,1)$ by (\ref{WQ2}) as $N\to\infty$.\\
For $I_2$ and for fixed $N>0$, by the dominated convergence theorem,
(\ref{EU6}) and the approximation of the identity (cf. \cite[p.23, (16)]{St0}), we have
\begin{align*}
I_2\leq\int^1_0\!\!\!\int_{B_N}\sup_{u\in B_R}|f_\eps(t,x,u)-f_{\eps'}(t,x,u)|^p\dif x\dif t\to 0,\ \ \eps,\eps'\to 0.
\end{align*}
Combining the above calculations and letting $\eps,\eps'\downarrow 0$ for (\ref{EEQ1}), we obtain
$$
\varlimsup_{\eps,\eps'\downarrow 0}\sup_{s\in[0,1]}\|w_{\eps,\eps'}(s)\|^p_{1-\frac{1}{p},p}=0,\ \
\varlimsup_{\eps,\eps'\downarrow 0}\int^1_0\|\nabla w_{\eps,\eps'}(s)\|^p_p\dif s=0.
$$
Hence, there exists a $u\in C([0,1];\mW^{1-\frac{1}{p},p}(\mR^d;\mR^m))\cap
L^p([0,1];\mW^{1,p}(\mR^d;\mR^m))$ such that
$$
\lim_{\eps\downarrow 0}\sup_{s\in[0,1]}\|u_\eps(s)-u(s)\|^p_{1-\frac{1}{p},p}
=0,\ \ \lim_{\eps\downarrow 0}\int^1_0\|\nabla u_\eps(s)-\nabla u(s)\|^p_p\dif s=0.
$$
Taking limits in $L^p$-space for equation (\ref{EY66}), it is easy to see that $u$ solves equation (\ref{Eq03}).
\end{proof}

\br
In this remark, we explain how to use the above results to the critical Hamilton-Jacobi equation
(cf. \cite{Dr-Im,Si3}). Let
$$
H(t,x,u,q):[0,1]\times\mR^d\times\mR^m\times\mM_{m\times d}\to\mR^m
$$
be a measurable and smooth function in $x,u,q$,
where $\mM_{m\times d}$ denotes the set of all real valued $m\times d$-matrices.
Consider the following Hamilton-Jacobi equation
\begin{align}
\p_t u=\cL u+H(t,x,u,\nabla u),\ \ u(0)=\varphi.\label{Ham}
\end{align}
Formally, taking the gradient we obtain
$$
\p_t \nabla u=\cL\nabla u+\nabla_x H(t,x,u,\nabla u)+\nabla_u H(t,x,u,\nabla u)\cdot\nabla u
+\nabla_q H(t,x,u,\nabla u)\cdot\nabla \nabla u.
$$
If we let
$$
w(t,x):=(u(t,x),\nabla u(t,x))^{\mathrm{t}},
$$
then
$$
\p_t w=\cL w+b(w)\cdot\nabla w+f(w),\ \ w(0)=(\varphi,\nabla\varphi)^{\mathrm{t}},
$$
where for $w=(u,q)$,
$$
b(t,x,w):=(0,\nabla_qH(t,x,u,q))
$$
and
$$
f(t,x,w):=(H(t,x,u,q),\nabla_x H(t,x,u,q)+\nabla_u H(t,x,u,q)\cdot q)^{\mathrm{t}}.
$$
Thus, we can use Theorems \ref{Th44} and \ref{Th4} to uniquely solve equation (\ref{Ham}) under
some assumptions on $H$ and $\varphi$.
\er

\end{document}